\documentclass[12pt]{article}
\pdfoutput=1
\usepackage{amsmath}
\usepackage{amssymb}
\usepackage{amsfonts}
\usepackage{amsthm}
\usepackage{amscd}
\usepackage{mathrsfs}
\usepackage{mathtools}
\usepackage{bm}
\usepackage{enumerate}
\usepackage{outlines} 
\usepackage{tabu} 
\usepackage{booktabs}
\usepackage{graphicx} 
\usepackage{grffile} 
\usepackage[margin=1.1cm, font=footnotesize]{caption} 
\usepackage{float} 
\usepackage{tikz} 
\usepackage{color}
\usepackage{xcolor}
\usepackage[square,sort,comma,numbers]{natbib}
\usepackage{algorithm}
\usepackage[noend]{algpseudocode}

\usepackage[colorlinks]{hyperref}
\hypersetup{final}
\usepackage{nomencl} 
\makenomenclature

\usepackage{subfigure}
\usepackage{hyperref}
\usepackage[nameinlink,capitalize]{cleveref}
\usepackage{authblk}

\newcommand{\rd}{\mathrm{d}}

\newcommand{\bI}{\mathbf{I}}

\newcommand{\cC}{\mathcal{C}}

\newcommand{\cG}{\mathcal{G}}
\newcommand{\cH}{\mathcal{H}}

\newcommand{\cK}{\mathcal{K}}
\newcommand{\cL}{\mathcal{L}}

\newcommand{\cN}{\mathcal{N}}

\newcommand{\cU}{\mathcal{U}}

\newcommand{\cW}{\mathcal{W}}

\newcommand{\bbC}{\mathbb{C}}

\newcommand{\bbR}{\mathbb{R}}

\newcommand{\bbZ}{\mathbb{Z}}

\newcommand{\hath}{\widehat{h}}

\newcommand{\tilA}{\widetilde{A}}

\newcommand{\tilK}{\widetilde{K}}

\newcommand{\tilU}{\widetilde{U}}

\newcommand{\tilW}{\widetilde{W}}

\newtheorem{theorem}{Theorem}

\newcommand{\p}{\partial}

\newcommand{\imag}{\mathrm{Im}}

\DeclarePairedDelimiter{\ceil}{\lceil}{\rceil}

\newcommand{\Ai}{A_{\text{in}}}
\newcommand{\Ao}{A_{\text{out}}}
\newcommand{\tilAi}{\tilA_{\text{in}}}
\newcommand{\tilAo}{\tilA_{\text{out}}}

\providecommand{\keywords}[1]
{
  \small	
  \textbf{\textit{Keywords---}} #1
} 

\title{Inverse Problem of Nonlinear Schr\"odinger Equation as Learning of Convolutional Neural Network}

\author[1]{Yiran Wang}
\author[2]{Zhen Li\thanks{Corresponding Author: 
lishen03@gmail.com}}
\affil[1]{Department of Mathematics, The Chinese University of Hong Kong, Hong Kong.}
\affil[2]{Theory Lab, Huawei Technologies Co., Ltd., Shenzhen, China.}

\begin{document}
\maketitle
\graphicspath{{pics/}}

\begin{abstract}
    In this work, we use an explainable convolutional neural network (NLS-Net) to solve an inverse problem of the nonlinear Schr\"odinger equation, which is widely used in fiber-optic communications. The landscape and minimizers of the non-convex loss function of the learning problem are studied empirically. It provides a guidance for choosing hyper-parameters of the method. The estimation error of the optimal solution is discussed in terms of expressive power of the NLS-Net and data. Besides, we compare the performance of several training algorithms that are popular in deep learning. It is shown that one can obtain a relatively accurate estimate of the considered parameters using the proposed method. The study provides a natural framework of solving inverse problems of nonlinear partial differential equations with deep learning.
\end{abstract}


\section{Introduction}

The nonlinear Schr\"odinger equation (NLSE) is a very special one among all nonlinear partial differential equations (PDE) in the sense that it has origin in many different phenomena of nature, such as propagation of lights in optical fibers and waveguides \cite{agrawal_nonlinear_2019}, Bose-Einstein condensates \cite{pitajevskij2003bose} and gravity waves on the surface of water \cite{Zakharov1968StabilityOP}, and it is also a subject that has attracted many mathematicians \cite{babelon_bernard_talon_2003,Bao2013}. Although not been aware of, the physics described by NLSE influences everyone's daily life, since optical networks are backbones underlying the modern internet \cite{agrawal_fiber-optic_2010}.

The work of this paper is initially motivated by applications of the NLSE in fiber-optic communications. This equation describes the propagation of light in optical fibers, which is used as the carrier of information. The solution of the inverse problem of this equation has potential applications such as optical fiber sensing, channel estimation and signal recovery in communication systems \cite{agrawal_fiber-optic_2010}.

\begin{figure}[H]
    \centering
    \includegraphics[width=0.5\textwidth]{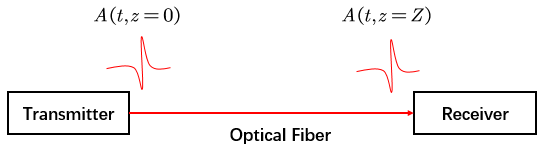}
    \caption{Schematic illustration of a fiber-optic communication system}
    \label{fig:fiber_link}
\end{figure}

The NLSE studied in this paper is
\begin{align}
\label{eq:NLSE}
    \frac{\p A}{\p z} =-\frac{i\beta}{2}\frac{\p^{2} A}{\p t^{2}} +i\gamma |A|^{2} A,
\end{align}
where $A:\bbR\times[0,Z]$ is the complex-valued amplitude of light field. The equation is viewed as an evolutionary equation in the spatial variable $z$. See \cref{fig:fiber_link} for an schematic illustration. A signal is injected from $z=0$ by a transmitter. It propagates through the fiber of length $Z$ and is received at $z=Z$ by a receiver.
While propagating in the optical fiber, the signal suffers from chromatic dispersion and Kerr nonlinear effect, which are described by the first and second term on the right hand side (RHS) of \cref{eq:NLSE},  respectively.
The readers are referred to \cite{agrawal_nonlinear_2019} for a comprehensive introduction of physical aspects of this equation.

The \textit{inverse problem} for the NLSE \eqref{eq:NLSE} considered in this paper is as follows. Suppose that we are given some data, which are samplings of the signal from both the transmitter and the receiver, the task is to determine the coefficients $\beta$ and $\gamma$ from the given data.

Our approach to this problem is differentiable programming, a version of explainable deep learning \cite{lecun_deep_2015,barredo_arrieta_explainable_2020,baydin_automatic_2018}. Artificial intelligence (AI) and artificial neural network have a long history \cite{mcculloch_logical_1943}. In the past ten years, it came back as deep learning and achieved great success in applications such as computer vision \cite{Deng_2009,gatys2016image,he_deep_2016}, machine translation \cite{Bahdanau2015NeuralMT} and game playing  \cite{Silver1140}. Recently, it is also developed as a tool for scientific research such as solving high-dimensional partial differential equations (PDE) in financial engineering, computing ground states of many-body quantum systems and predicting folding structures of proteins \cite{han_solving_2019,choo_fermionic_2020,hermann_deep-neural-network_2020,pfau_ab_2020}. Particularly, they are also used to solve inverse problems \cite{li2020nett,adler2017solving,jin2017deep,aggarwal2018modl,adler2017learning,mukherjee2021adversarially,lunz2021learned,herzberg2021graph}.

However, deep learning is criticized by some serious researchers as `alchemy' because of its lack of theoretical understanding \cite{hutsonmay_3_ai_2018}, i.e., it is generally unexplainable. This difficulty has been considered by more and more people as a bottleneck for systematical improvement of AI technology, and also an obstacle for its applications in scientific research and in safety- and security- related applications such as self-driving, AI-aided healthcare and communication networks. The issue of explainability is also relevant to our task, i.e., solving inverse problems.

Viewed in the framework of machine learning, the unexplainability of deep learning lies in three major aspects. The first one is the representation (expressive) power of model. It is unknown whether the neural network model can approximate the right function class efficiently as required by a specific task. The second aspect is training. Even if the function class of the task can be approximated well by the model in theory, it is unknown whether the right parameters can be practically found by optimizing the loss function, because the optimization problem is of high-dimensional and non-convex. Finally, the generalization ability of deep neural network is unpredictable and unexplainable, i.e., it is not known when and why a trained model performs well (or bad) on test data. 

The representation problem has been studied in theory for a long time for classical function classes such as Sobolev spaces \cite{cybenko_approximation_1989,hornik_approximation_1991,mhaskar_neural_1996, yarotsky18a,yarotsky2019phase}. Error bounds are obtained for ReLU deep neural networks. But there is still no systematic way to design good architecture for a specific task.
Training and generalization theory are still largely not exist. Only a few results appear in recent years. For example, \cite{allen-zhu19a,allen-zhu19b} prove that the stochastic gradient descent (SGD) can achieve global optimum in the training of over-parameterized neural networks. Alternative to develop theory for general deep neural networks, one can integrate prior knowledge of a specific task into the design of neural network from the beginning. Such models are expected to have inherent explainability. Differentiable programming is a straightforward way to realize this idea.

In this paper, we construct a highly explainable neural network model from NLSE and formulate the inverse problem as a machine learning problem of the model. The formulation follows the procedure of differentiable programming. It is guaranteed by mathematical theory that the neural network with sufficiently large depth and width can approximate the solution of the problem arbitrarily well, and the error bound can be obtained explicitly. 
Then we visualize the loss landscape of the learning problem and observed that there is an unique global minimizer. We discuss the effects of model architecture and data on the estimation error of this optimal solution. 
Finally, we compare performance of popular training algorithms on the inverse problem of NLSE. We assume that we have sufficient prior knowledge about the solution, i.e., good initial guess, hence the problem is approximately convex and can be solved in a straightforward way. The general non-convex optimization problem is left as future work. The comparison are carried out with different settings of noise: noiseless data, noisy data and denoised data. Notice that we focus on the inverse problem in this paper and will not discuss the generalization problem, which is the goal of machine learning. That will be our future work.

\textbf{Related works.}
Recently, deep learning attracts lots of attention from mathematicians. In particular, many efforts have been made to solve inverse problems using neural networks \cite{li2020nett,adler2017solving,jin2017deep,aggarwal2018modl}. So far as we know, in all of these works, neural networks are used as an approximate mapping from data to the desired quantities. Our approach is fundamentally different. We transform the NLSE itself into a neural network, and the inverse problem becomes a learning problem.

Besides, some researchers have tried to relate DNN with PDE. Especially, there are several works on explaining neural network by PDE \cite{e_proposal_2017, LIZhen2017,long_pde-net_2018,ruthotto2019deep,wang2020mesh,kutyniok2019theoretical,elbrachter2018dnn}. So far as we know, they all use the similarity between residual network (ResNet) and Euler finite difference as a bridge to connect DNN and PDE \cite{he_deep_2016}. The current work provides a natural way to interpret PDE as a continuous version of non-residual DNN. This connection might be useful for the theory of deep learning.

\textbf{Structure of the paper.}
In \cref{ss:formulation}, the forward and inverse problems of the NLSE are described. Then the explainable neural network is constructed and the inverse problem is reformulated as a learning problem of the neural network. In \cref{ss:landscape}, the effects of model architecture and data on the solution is empirically studied. In \cref{ss:solving}, the performance of several training algorithms is compared. Finally, we conclude the paper in \cref{ss:discussion}.

\section{Formulation of the Problem} 
\label{ss:formulation}

\subsection{Abstract Description}
\label{ss:abs_formulation}
The main task of this paper is to solve an inverse problem of the NLSE \cref{eq:NLSE}, i.e., to estimate the coefficients of the equation from given data. To solve the inverse problem, one needs to solve the forward problem first, which is the initial value problem (IVP, Cauchy problem) of \cref{eq:NLSE} with initial condition at $z=0$:
\begin{align}
\begin{cases}
\label{eq:NLSE_IVP}
    \frac{\p A}{\p z} =-\frac{i\beta}{2}\frac{\p^{2} A}{\p t^{2}} +i\gamma |A|^{2} A, & t\in\bbR, z\in [0,Z]\\
    A(t,0) = \Ai(t), & t\in\bbR,
\end{cases}
\end{align}
where $\Ai:\bbR\to\bbC$. The signal propagates according to \eqref{eq:NLSE} and is received at $z=Z$ as $\Ao:\bbR\to\bbC$, which is defined by
\begin{align}
    \Ao(t) = A(t,Z), \quad t\in\bbR.
\end{align}
The observed data is the pair $(\Ai,\Ao)$.
The well-posedness of the IVP \cref{eq:NLSE_IVP} is established under various conditions. According to these theoretical results, we can take the light field $A\in H^{s}_{t}C_{z}(\bbR\times[0,Z]\to\bbC)$ and $\Ai, \Ao\in H^{s}(\bbR\to\bbC)$ for any $s>1/2$, where $H^{s}$ is the standard notation for Sobolev space. The notation $H^{s}_{t}C_{z}$ means that $A$ is $H^{s}$ in the variable $t$ and continuous in the variable $z$.
For technical reasons, we assume that $s\ge 5$.
For readers who are interested in the details, we recommend \cite{Bao2013,tao_nonlinear_2006}. 

The inverse problem is to find the coefficients $\beta$ and $\gamma$ in \cref{eq:NLSE} using the observed data. Generally, these coefficients are themselves functions of $(t,z)$. Here we put emphasis on the method and restrict ourselves to the simplest case that $\beta$ and $\gamma$ are constants, i.e., they do not depend on $(t,z)$. As long as the method for the simplest case is established, it is possible to extend the method to more general settings. 

Suppose that the exact propagation operator of the NLSE \cref{eq:NLSE} is $\cU(z;\beta,\gamma)$ such that for any $z\in\bbR$ there is $A(\cdot,z)=\cU(z;\beta,\gamma)A(\cdot,0)$. Then our inverse problem can be formulated as an optimization problem
\begin{align}\label{eq:opt}
    \min_{\beta,\gamma} L(\cU(Z;\beta,\gamma)\Ai, \Ao),
\end{align}
where $L$ is an error metric. It is not necessarily a metric in the sense of distance. For example, a common choice of $L$ is the $L^{2}$-distance $L(f,g) = \|f-g\|_{L^{2}}$. 

In the setting of numerical computation, only discrete and finite samplings of $\Ai$ and $\Ao$ can be collected. Let's consider the signals on $[0,T]$. Given a set of sampling time $\{t_{n}=n \tau\}_{n=0}^{N-1}$, our data is $\{(\Ai(t_{n}), \Ao(t_{n}))\}_{n=0}^{N-1}$. Here $N$ is the number of samples, and $\tau=T/N$ is the sampling period. For convenience of later use, we denote $\tilAi[n]=\Ai(t_{n})$ and $\tilAi=(\tilAi[n])_{n=0}^{N-1}$. Similarly for $\tilAo$. In practical scenarios, noise is inevitable. To make the problem setting more realistic, we can add noise to $\tilAi, \tilAo$. Their are different methods to add noise, which will be specified in \cref{ss:landscape,ss:solving}.
Generally, $\cU$ has no analytical expression and it is necessary to solve the IVP \cref{eq:NLSE_IVP} numerically. Let $\tilU$ be a numerical approximation of $\cU$, then an estimate of the true parameters $(\beta^{\dag},\gamma^{\dag})$ is a solution of the optimization problem
\begin{align}\label{eq:opt_comp}
    \min_{\beta,\gamma} L(\tilU\tilAi, \tilAo),
\end{align}
denoted as $(\beta^{*},\gamma^{*})$. 

Considering the accuracy of solution, there are two sources of error. The one is related to the model $\tilU$ itself. Another is related to the data $(\tilAi, \tilAo)$. As a result, a global minimizer of \cref{eq:opt_comp}, is not necessarily equal to the ground truth $(\beta^{\dag},\gamma^{\dag})$ that has been used to generate data. Later we will come back to analyse these errors.

\subsection{Convolutional Neural Network Model}\label{ss:CNN}
Convolutional neural network (CNN) is a type of neural network whose linear operators in (some of) its layers are convolution operators. CNN is one of the most important achievements in the development of AI algorithms \cite{fukushima_neocognitron_1980,waibel_phoneme_1989,lecun_backpropagation_1989}. It is widely used in application areas such as computer vision and speech processing \cite{lecun_deep_2015,he_deep_2016,abdel-hamid_convolutional_2014}. It is also applied to inverse problems in imaging \cite{jin_deep_2017}. In this paper, we use CNN to solve the inverse problem of PDE. More specifically, we use a CNN as an approximate propagation operator $\tilU$ of the NLSE \cref{eq:NLSE}. 

Our CNN is derived from a popular numerical method for solving IVP of NLSE: the split-step Fourier method (SSFM) \cite{hardin_application_1973}. It is an application of the idea of operator splitting \cite{glowinski_splitting_2016,strang_1968}. Here we only describe the basic idea and major results of the method and omit the details of the derivation. The core of SSFM is the operator splitting method.
The RHS of the NLSE \cref{eq:NLSE} can be split into two parts: the linear part and the nonlinear part. Written in abstract form,
\begin{align}
    \frac{\p A}{\p z} = \cL A + \cN A,
\end{align}
where the linear operator $\cL$ and the nonlinear operator $\cN$ are defined by
\begin{align}
    \cL A &:= -\frac{i\beta}{2}\frac{\p^{2} A}{\p t^{2}}, \\
    \cN A &:= i\gamma |A|^{2} A.\label{operator}
\end{align}
The linear and nonlinear parts of the NLSE can be solved separately. The propagation operator of the linear part is given by
\begin{align}
    A(t, z+\zeta) = \int_{s\in\bbR} K(s;\beta\zeta) A(t-s, z) \rd s, \quad \forall t \in\bbR,
\end{align}
where the kernel $K(\cdot;\eta)$, parameterized by $\eta\in\bbR$, is defined by
\begin{align}
    K(t;\eta):= \sqrt{\frac{i}{2\pi \eta}}
        \exp\left(-\frac{i t^{2}}{2\eta} \right), \quad\forall t\in\bbR.
\end{align}
Formally, we denote
\begin{align}\label{eq:linear_solver}
    A(\cdot, z+\zeta) = \cK(\beta\zeta) A(\cdot, z) := K(\cdot;\beta\zeta) * A(\cdot, z),
\end{align}
where $*$ is the convolution with respect to $t$.
The propagation operator of the nonlinear part is given by
\begin{align}
    A(t, z+\zeta) = \kappa(A(t,z);\gamma\zeta),\quad\forall t\in\bbR.
\end{align}
where the nonlinear function $\kappa(\cdot;\eta)$, parameterized by $\eta\in\bbR$, is defined by
\begin{align}
    \kappa(w;\eta) := w \exp \left(i \eta |w|^{2}\right), \quad\forall w\in\bbC.
\end{align}
In the above expressions, $\zeta\in\bbR$ is the distance of the propagation.

To solve the IVP \cref{eq:NLSE_IVP}, the spatial interval $[0,Z]$ is partitioned as $0=z_{0}<z_{1}<\dots<z_{M}=Z$. For simplicity, let the partition be uniform, i.e., for any $m\in[1:M]$ there is $z_{m}-z_{m-1}=\zeta>0$. Here $[1:M]$ is a short notation for the set of integers $\{m : 1\le m \le M\}$. The propagation operator $\cU$ for the IVP \cref{eq:NLSE_IVP} is approximated by a stack of alternating compositions of the convolution operator $\cK$ and pointwise nonlinear operator $\kappa$. There are different choices for arrangement of the operators. We adopt the Strang splitting \cite{strang_1968}
\begin{align}
    \cU(Z; \beta, \gamma) \approx U(\beta, \gamma; M, \zeta) := \left( \cK\left(\frac{\beta\zeta}{2}\right) \circ \kappa(\cdot;\gamma\zeta) \circ \cK\left(\frac{\beta\zeta}{2}\right) \right)^{\circ M}.
\end{align}
Here $\circ$ is the composition of operators and $(\cdot)^{\circ M}$ means composition with itself for $M$ times. The operator $U(\beta, \gamma; M, \zeta)$ is a semi-discretization of $\cU(Z; \beta, \gamma)$.

We are left to discretize $\cK$ in order to solve \cref{eq:NLSE_IVP} numerically. We adopt the trapezoidal quadrature rule. To be consistent with our earlier notations of data, for any $m\in[0:M]$, let $\tilA_{m}[n]$ be the value of $A$ at $(t_{n}, z_{m})$ and denote $\tilA_{m}=(\tilA_{m}[n])_{n=0}^{N-1}$. 
The discrete convolution kernel $\tilK(\eta)$ is defined by
\begin{align} \label{eq:discrete_kernel}
    \tilK(\eta,\tau)[k]:=\sqrt{\frac{i}{2\pi \eta}}
        \exp\left(-\frac{i (k\tau)^{2}}{2\eta} \right), \quad\forall k\in \left[-\ceil*{\frac{N}{2}}:\ceil*{\frac{N}{2}}\right].
\end{align}
Denote the corresponding linear operator by $\widetilde{\cK}(\eta,\tau)$, which is a discretization of $\cK(\eta)$. Then the convolution \cref{eq:linear_solver} can be approximated by
\begin{align}
    \cK(A(\cdot, z_{m});\beta\zeta)(t_{n}) 
    \approx & (\widetilde{\cK}(\beta\zeta,\tau)\tilA_{m})[n] \\
    = & (\tilK(\beta\zeta,\tau) * \tilA_{m})[n] \\
    = & \sum_{l=0}^{N-1} \tilK(\beta\zeta,\tau)[n-k] \tilA_{m}[l],
\end{align}
for any $m\in [0:M-1]$ and $n\in[0:N-1]$. Here we abuse $*$ as the notation of discrete convolution. In the summation, the terms which are not defined for some $k$'s are understood as zeros.
Then we have
\begin{align}
\label{eq:splitting}
    \tilA_{m+1}[n] := \left(\left( \widetilde{\cK}\left(\frac{\beta\zeta}{2},\tau\right) \circ \kappa(\cdot;\gamma\zeta) \circ \widetilde{\cK}\left(\frac{\beta\zeta}{2},\tau\right) \right) \tilA_{m}\right)[n]
\end{align}
This recursive formula gives us the numerical algorithm of SSFM for the IVP \cref{eq:NLSE_IVP}. The operator
\begin{align}
    \tilU(\beta, \gamma; M, N, \zeta, \tau)
    := \left( \widetilde{\cK}\left(\frac{\beta\zeta}{2},\tau\right) \circ \kappa(\cdot;\gamma\zeta) \circ \widetilde{\cK}\left(\frac{\beta\zeta}{2},\tau\right) \right)^{\circ M}
    \label{eq:SSFM}
\end{align}
is a discrete approximation of the propagation operator $\cU(Z; \beta, \gamma)$ for the NLSE \cref{eq:NLSE}.

\begin{figure}[H]
    \centering
    \includegraphics[width=0.8\textwidth]{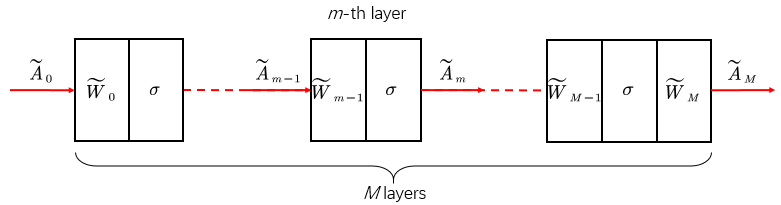}
    \caption{Schematic illustration of the architecture of NLS-Net.}
    \label{fig:SSFM}
\end{figure}

The expression \cref{eq:SSFM} is a stack of alternating compositions of the convolution operators and pointwise nonlinear operators. It can be equivalently written as
\begin{align}
    \tilU(\beta, \gamma; M, N, \zeta, \tau)
    = \left( \frac{1}{\sqrt{\gamma\zeta}}\widetilde{\cK}\left(\frac{\beta\zeta}{2},\tau\right) \circ \sigma \circ \sqrt{\gamma\zeta} \widetilde{\cK}\left(\frac{\beta\zeta}{2},\tau\right) \right)^{\circ M},
\end{align}
where $\sigma(w):=\kappa(w;1)$. This is exactly in the form of a $M$-layer complex-valued CNN with the activation function $\sigma$:
\begin{align}
    \widetilde{\cG}(\Theta)
    = \widetilde{\cW}_{M} \circ
    \sigma \circ \widetilde{\cW}_{M-1}
    \circ\dots\circ
    \widetilde{\cW}_{1} \circ
    \sigma \circ \widetilde{\cW}_{0},
    \label{eq:NLS-Net}
\end{align}
where $\Theta:=\{\tilW_{m}\}_{m=0}^{M}$ is the collection of convolution kernels corresponding to the linear operators $\{\widetilde{\cW}_{m}\}_{m=0}^{M}$. The $M$ and $N$ are usually called the depth and width of the neural network, respectively. In this paper, the sizes of the kernels are all set to $N$, hence the width of each layer is $N$. It is straightforward to generalize the architecture in this paper to neural networks with different width in each layer. To do that, one only needs to use different time partitions for $A(t, z_{m})$ for different $m$.  
For convenience, we call this CNN as \textit{nonlinear Schr\"odinger network (NLS-Net)}.
See \cref{fig:SSFM} for an illustration of NLS-Net.

The NLS-Net is derived from SSFM, whose convergence properties are proved rigorously \cite{Bao2013}. In contrast to most of neural networks, our model is a deductive rather than heuristic approximation of the propagation operator $\cU$, hence it is inherently explainable. It means that the NLS-Net indeed express the propagation operator of the NLSE \cref{eq:NLSE} if the parameters are set properly. This can be stated formally as a theorem \cite{Bao2013}:
\begin{theorem}[Expressive Power] \label{thm:convergence_rate}
Let the input data $\tilAi=(A(t_{n},0))_{n=0}^{N-1}$ and the target data $\tilAo=(A(t_{n},Z))_{n=0}^{N-1}$ be given as a sampling of the initial and terminal value of the IVP \cref{eq:NLSE_IVP}. Suppose that the initial value $A(\cdot,0)\in H^{s}(\bbR\to\bbC)$ for $s\ge 5$. Then there exist parameters $\Theta$ of the neural network $\widetilde{\cG}(\Theta)$, which is defined by \cref{eq:NLS-Net}, such that
\begin{align}
    \|\widetilde{\cG}(\Theta)\tilAi-\tilAo \| \lesssim \zeta^{2} + \tau^{s}.
\end{align}
\end{theorem}
\begin{proof}
In the NLS-Net defined by \cref{eq:NLS-Net}, let the convolution kernels be given by
\begin{align}\label{eq:reparameterization}
\begin{cases}
    \tilW_{0} = \sqrt{\gamma\zeta} \tilK\left(\frac{\beta\zeta}{2},\tau\right),\\
    \tilW_{1} = \dots = \tilW_{M-1} = \tilK\left(\beta\zeta,\tau\right),\\
    \tilW_{M} = \frac{1}{\sqrt{\gamma\zeta}}\tilK\left(\frac{\beta\zeta}{2},\tau\right).
\end{cases}
\end{align}
Then the NLS-Net becomes the discrete propagation operator \cref{eq:SSFM} provided by SSFM. According to the convergence results of SSFM \cite{Bao2013}, the approximation error of the terminal value computed by SSFM is bounded by $\zeta^{2} + \tau^{s}$ from above under conditions of the theorem.
\end{proof}
Notice that $M=Z/\zeta$ and $N=T/\tau$, thus the \cref{thm:convergence_rate} characterizes the expressive power of NLS-Net by giving an explicit upper bound of the approximation error of NLS-Net in terms of its depth and width for fixed $T$ and $Z$.
It is a foundation for using the NLS-Net \cref{eq:NLS-Net} to solve the inverse problem of the NLSE \cref{eq:NLSE}. If the parameters are set equal to the ground truth that are used in the NLSE \cref{eq:NLSE}, from which the data is generated, then the output of the neural network will agree with the target data up to numerical error of SSFM, which is controlled by data and the size of NLS-Net.

However, this result is unsatisfactory in the view of machine learning, because the upper bound $\zeta^{2} + \tau^{s}$ does not converge to zero for fixed depth or fixed width. Besides, it gives no information about the lower bound of approximation error. These issues will be touched later in an empirical investigation of loss landscape in \cref{ss:landscape}.

\subsection{Supervised Learning of NLS-Net}
In \cref{ss:abs_formulation}, we have shown that the inverse problem for NLSE can be solved through the optimization problem \cref{eq:opt_comp}. In \cref{ss:CNN}, we have designed a NLS-Net \cref{eq:NLS-Net} to approximate $\cU(Z; \beta,\gamma)$. 
Therefore, given the input data $\tilAi$ and target data $\tilAo$, our task is to train the NLS-Net model to minimize the gap between the output $\widetilde{\cG}(\Theta)\tilAi$ and the target $\tilAo$. This is a supervised learning problem for the NLS-Net.

\begin{figure}[!htbp]
    \centering
    \includegraphics[width=\textwidth]{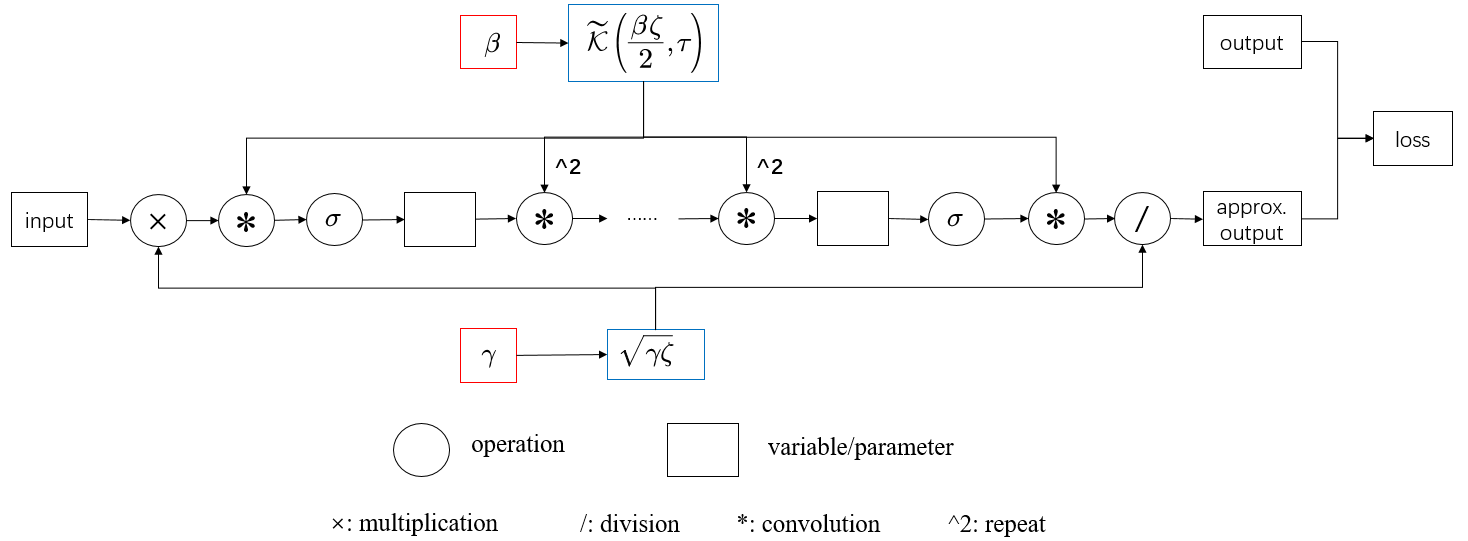}
    \caption{Forward Computation graph of NLS-Net. The squares stand for variables and parameters. The circles stand for operations. The directions of arrows indicate the inputs and outputs of operations.}
    \label{fig:computation_graph}
\end{figure}

In principle, all the linear operators in the NLS-Net \cref{eq:NLS-Net} are completely trainable. 
But in the particular scenario of solving the inverse problem of \cref{eq:NLSE}, we only need to train the parameters $\beta$ and $\gamma$. Therefore, we re-parameterize the convolution kernels by $\Theta = \theta(\beta, \gamma; \zeta, \tau)$, which is defined by \cref{eq:reparameterization}, and freeze all the parameters in the kernels, except $\beta$ and $\gamma$. The re-parameterized NLS-Net model is defined by
\begin{align}\label{eq:reparameterized_model}
    \widetilde{\cH}(\beta, \gamma; \zeta, \tau):= & \widetilde{\cG}(\theta(\beta, \gamma; \zeta, \tau)) \\
    =& \frac{1}{\sqrt{\gamma\zeta}}\widetilde{\cK}\left(\frac{\beta\zeta}{2},\tau\right) \circ \sigma \circ \widetilde{\cK}(\beta\zeta,\tau) \circ \cdots
    \circ \widetilde{\cK}(\beta\zeta,\tau)
    \circ \sigma \circ \sqrt{\gamma\zeta} \widetilde{\cK}\left(\frac{\beta\zeta}{2},\tau\right).
\end{align}
The corresponding optimization problem \cref{eq:opt_comp} is reduced to
\begin{align}\label{eq:opt_reduced}
 \min_{\beta,\gamma} J(\beta, \gamma; \tilAi, \tilAo)
:= \min_{\beta,\gamma} L(\widetilde{\cH}(\beta, \gamma)\tilAi, \tilAo).
\end{align}
Here $J$ is the objective function, which is often called loss function. A global minimizer of $J$ is an optimal estimate of the true parameters $(\beta^{\dag},\gamma^{\dag})$, denoted as $(\beta^{*},\gamma^{*})$.
The structure of the model and the optimization problem can be visualized by the computation graph as shown in \cref{fig:computation_graph}.
Recall that the hyper-parameters $M, N$ are the depth and width of the NLS-Net, respectively. $\zeta=Z/M$ and $\tau=T/N$ are also known from data. Hence they are implicitly included in the definition of $J$. They are set before the optimization problem is solved. Only $\beta, \gamma$ are optimization variables. In the above and following expressions of $\widetilde{\cH}$, we drop the hyper-parameters $\zeta,\tau$ whenever it is clear from contexts.

The loss function $J$ defined above is to be specified further. It depends on the choice of the metric $L$ as well as the given data $(\tilAi, \tilAo)$. Here we choose the squared $L^{2}$-distance $L(f,g) = \|f-g\|_{L^{2}}^{2}$, but normalized with target function. Hence the discrete version of loss function is given by
\begin{align}\label{eq:full_loss}
    J(\beta, \gamma; \tilAi, \tilAo)
= \frac{\left\|\widetilde{\cH}(\beta, \gamma)\tilAi- \tilAo\right\|_{2}^{2}}{\|\tilAo\|_{2}^{2}}
=\frac{\sum_{n=0}^{N-1} \left|\left(\widetilde{\cH}(\beta, \gamma)\tilAi\right)[n] - \tilAo[n] \right|^{2}}{\sum_{n=0}^{N-1} |\tilAo[n]|^{2}}.
\end{align}

Notice that we only need one pair of data $(\tilAi,\tilAo)$ for the optimization of $J$, which is in the form of one-shot learning. One can also use a training dataset with multiple pairs of $(\tilAi,\tilAo)$ as in ordinary machine learning. However, different pairs of data should be similar to each other in the statistical sense, since our data are generated from random symbol sequence $\{a_{k}\}$. Therefore, using multiple short data should be effectively equivalent to using a single long data. The two approaches should have similar properties and performance in our task.

\subsection{Data}
\label{ss:data}

The NLSE \cref{eq:NLSE} considered in this work serves as a model for the transmission of signals in optical fibers. Hence we mainly consider the form of data in this scenario. The inverse problem of other forms of the NLSE, such as the Gross-Pitaevskii equation (GPE) describing the dynamics of Bose-Einstein condensates, can also be solved using our method. But they are out of the scope of the current work. All the data used in this paper are synthetic. We first describe the form of data analytically, then we describe the procedure for sampling and adding noise. 

\begin{figure}[ht]
    \centering
    \subfigure[$\hath(f)$ with different $\rho$]
      { \includegraphics[width=0.32\linewidth]{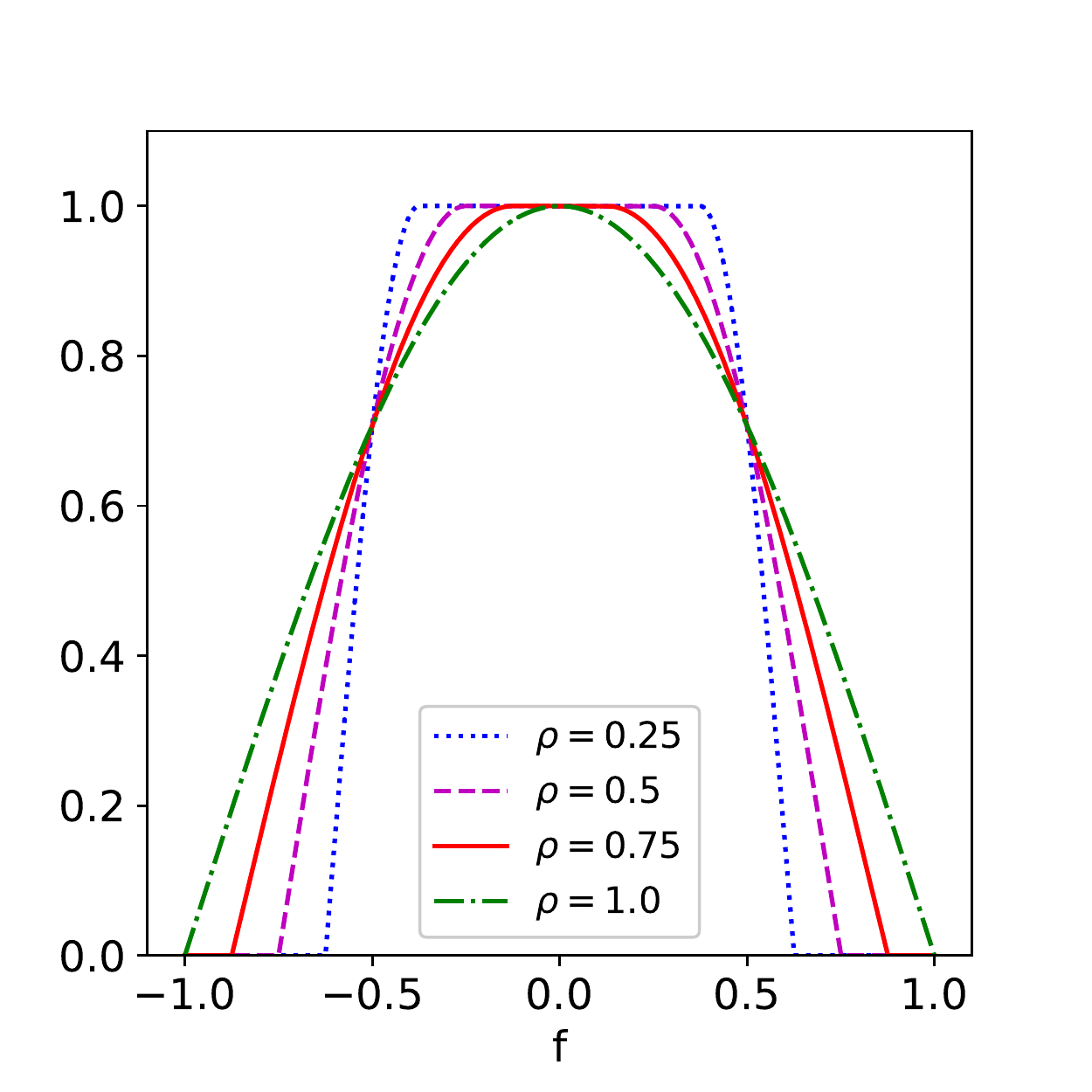}\label{fig:pulse_freq}}
     \subfigure[$h(t)$ with different $\rho$]
	{\includegraphics[width=0.32\linewidth]{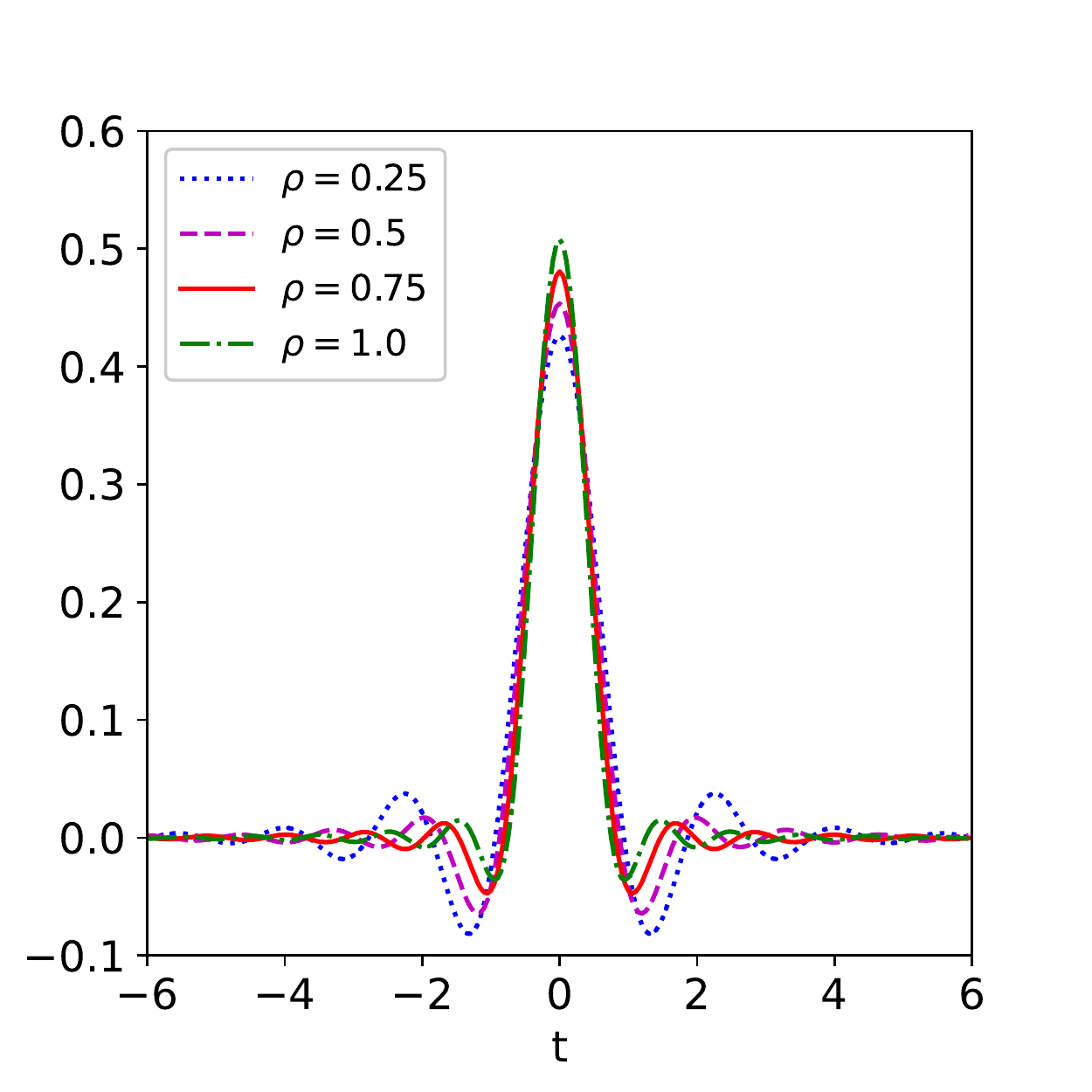}
	\label{fig:pulse_time}}
    \subfigure[16QAM constellation $\cC$]
	{\includegraphics[width=0.32\linewidth]{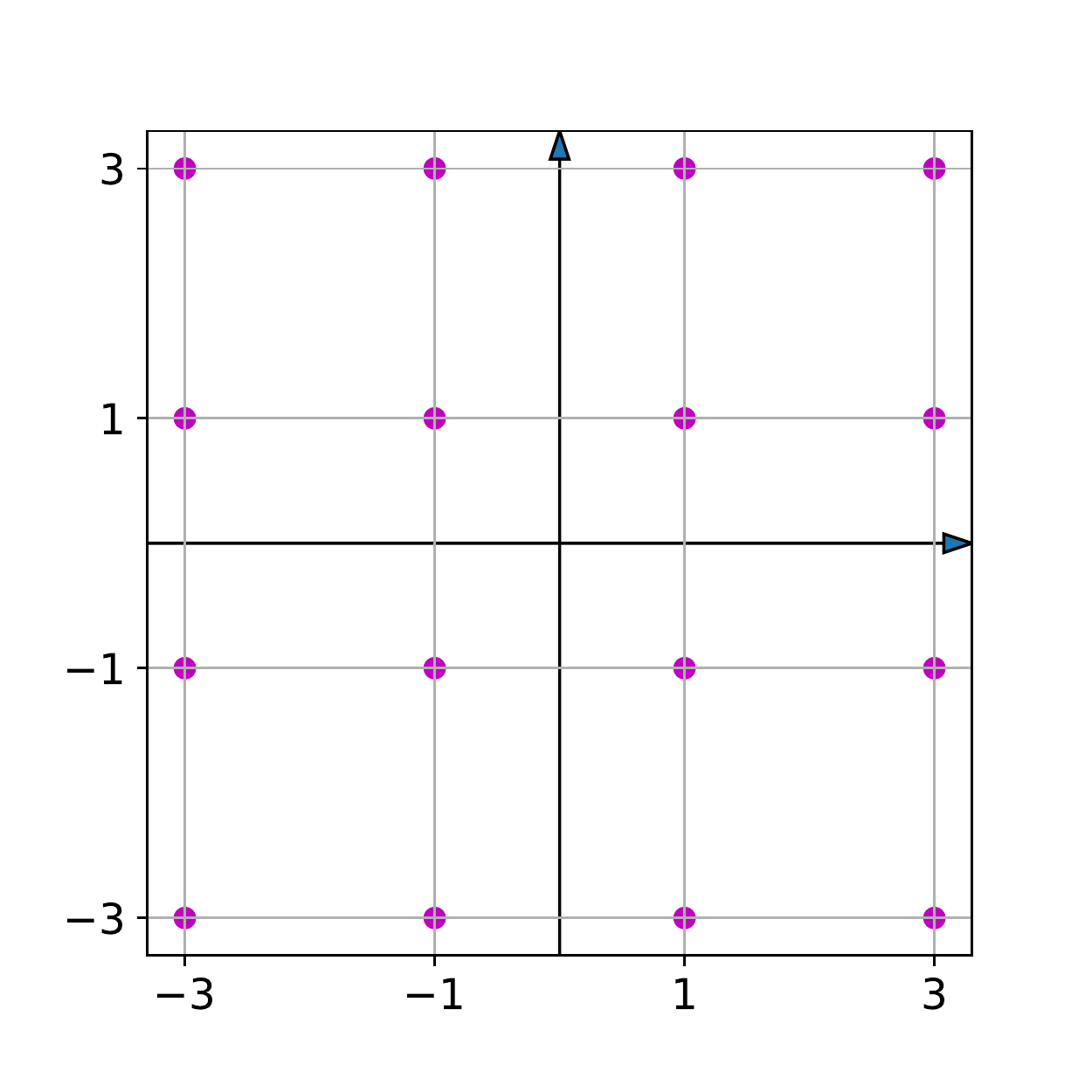}
	\label{fig:16qam}}
	\caption{An illustration of $h$ and $\cC$.}
	\label{fig:pulse_constellation}
\end{figure}

\begin{figure}[!ht]
    \centering
	\includegraphics[width=0.8\textwidth]{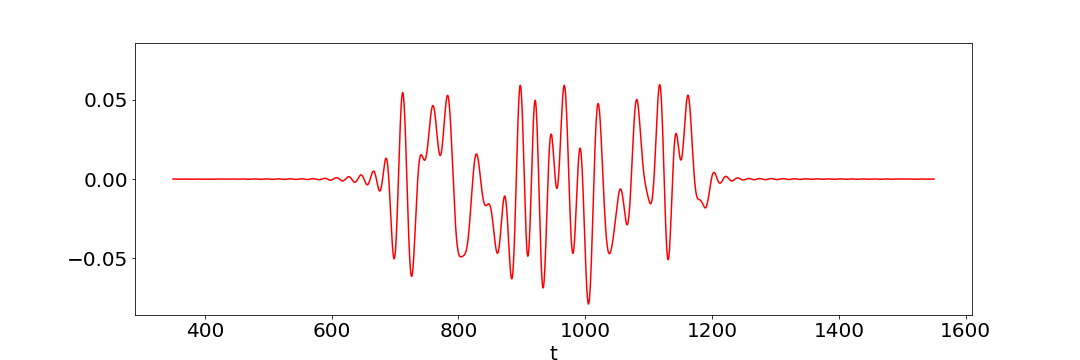}
	\caption{An example of 16QAM signal. Only the real part is shown. The imaginary part has similar style of waveform.}
    \label{fig:illustration_A_0}
\end{figure}

The initial data of the IVP \cref{eq:NLSE_IVP} is in the form of
\begin{align}\label{eq:Ai}
    \Ai(t) = \sqrt{P}\sum_{k\in\bbZ} a_{k} h(t-kT_{s}), \quad\forall t\in\bbR. 
\end{align}
Here $P>0$ is a scaling factor, which is set to 1 throughout this paper. The coefficients $a_{k}$'s, called symbols, are randomly chosen from a finite set $\cC\subset\bbC$, called constellation, with uniform probability. The input data $\Ai$ can be regarded as been generated from the sequence $\{a_{k}\}_{k\in\bbZ}$, called symbol sequence, which is a representation of information. The $T_{s}>0$ is called symbol period and its inverse is called symbol rate. For example, a signal $\Ai$ is called 100 GBaud if its symbol rate is 100 per picosecond. In the area of optical communication, there are different choices of the set $\cC$ and the pulse function $h:\bbR\to\bbR$. We follow the most popular choices, which are specified in the following. Let $h$ be a root-raised-cosine (RRC) filter, whose frequency-domain description is 
given by
\begin{align}\label{eq:pulse_shape}
\hath(f)=
\begin{cases}
1, & |f| \leq \frac{1-\rho}{2 T_{s}}, \\
\cos \left(\frac{\pi T_{s}}{2 \rho}\left(|f| - \frac{1-\rho}{2 T_{s}}\right)\right), & \frac{1-\rho}{2 T_{s}}<|f| \leq \frac{1+\rho}{2 T_{s}}, \\
0, & \text{otherwise},
\end{cases}
\end{align}
where $f$ is the frequency and $\rho$ is called roll-off factor. 
See \cref{fig:pulse_freq} and \cref{fig:pulse_time} for illustrations of $h$ in frequency and time domain respectively.
Let the set $\cC$ be defined by $\{\pm (2m+1) \pm i (2n+1)\}_{0\le m,n\le1}$, which is called 16QAM (quadrature-amplitude modulation) constellation because it contains 16 grid points. 
See \cref{fig:16qam} for an illustration of $\cC$.
One example of the initial signal $\Ai$ is illustrated in \cref{fig:illustration_A_0}.
Suppose that the ground truth of the parameters to be estimated is $(\beta^{\dag},\gamma^{\dag})$. Then we can solve the IVP \cref{eq:NLSE_IVP} with coefficients  $(\beta^{\dag},\gamma^{\dag})$. By definition, The received signal $\Ao$ is the solution of the IVP at position $z=Z$, that is $\Ao=A(\cdot,Z)$. 

In the computational setting, we can only represent finite symbols, i.e., there are only finite terms in the sum \cref{eq:Ai}. Hence the discrete version of data is a sampling of \cref{eq:Ai}, which is in the form of 
\begin{align}\label{eq:Ai_comp}
    \tilAi[n] = \sqrt{P}\sum_{k=1}^{N_{s}} a_{k} h(n\tau-kT_{s}) + \epsilon, \quad\forall n=0,\dots,N-1. 
\end{align}
Here $N_{s}$ is the number of symbols in the data, which is also referred to the length of data. It is convenient to set $\tau$ such that $T_{s}/\tau$ is an integer. It is the sampling rate measured in the unit of samples per symbol (sps). 
The $\epsilon$ is the noise added during sampling. The noise is drawn from $N$-dimensional white Gaussian distribution $\cN(\mathbf{0}, \delta^{2}\bI)$, where $\delta \ge 0$ will be specified in numerical experiments. A more convenient way of characterize the level of noise is to use the signal-to-noise-ratio (SNR). It is defined as the ratio of the $L^{2}$-norms of the noiseless signal and the noise. Noiseless data is considered as a special case of noisy data with infinite SNR.

The output data $\tilAo$ is generated by solving the IVP \cref{eq:NLSE_IVP} with noiseless initial data $\tilAi$ using SSFM. To avoid artificial boundary effect of SSFM, we pad zero symbols to $\tilAi$. The $a_{k}$'s near the boundary are set to zero. And $N_{s}$ stands for nonzero $a_{k}$'s. Alternatively, one can also set $a_{N_{s}}=a_{1}$ instead of using zero symbols.
The same level of noise is added to the output data $\tilAo$ in a similar way. In practice, there is also noise during the transmission, but this is omitted here.

\section{Loss Landscape}\label{ss:landscape}

There are two major aspects for a general optimization problem: loss landscape (landscape of the loss function) and optimization algorithm. For non-convex optimization problems as here, both aspects are crucial. 
In this section, we first investigate major characters of the landscape of $J$ for a typical choice of hyper-parameters. Then we study the effects of model architecture and data on the optimal estimate of the true parameters, i.e., the global minimizer of the loss landscape. The error of the optimal estimate will be discussed in terms of expressive power of NLS-Net and randomness of data. The data used in this section is noiseless by default.

\subsection{Major Characters of Loss Landscape} 
\label{ss:landscape_characters}

For training problem of a general DNN, there are two major difficulties: high dimensionality of optimization space and non-convexity of loss function, which make the problem both unexplainable in theory and difficult in practice. Fortunately, in our problem, the optimization space is low-dimensional, since most parameters in our NLS-Net are frozen and only $\beta,\gamma$ are regarded as optimization variables. This advantage make it possible to visualize and to have a much clearer and thorough investigation of the loss landscape.

\begin{figure}[!htbp]
    \centering
    \subfigure[$J$ in small region]
      { \includegraphics[width=0.4\textwidth]{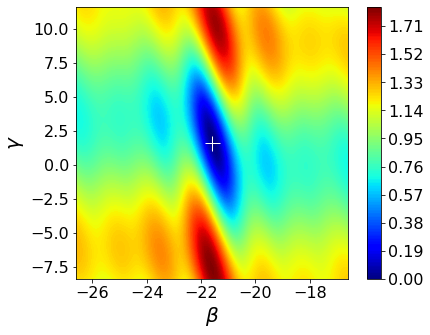}\label{l2-local}}
    \subfigure[$J$ in large region]
      { \includegraphics[width=0.4\textwidth]{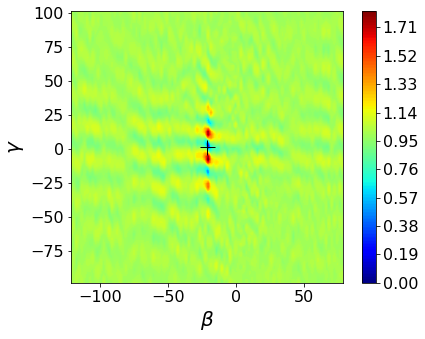}\label{l2-global}}
    \caption{Visualization of loss landscape. Colors indicate values of $J$. The location of ground truth is annotated by ``+'' in each panel.}
    \label{fig:loss_landscape}
\end{figure}

We compute the value of $J$ on a set of grids of $(\beta, \gamma)$ and visualize the loss landscape in \cref{fig:loss_landscape}. The major hyper-parameters for generating the figures are listed in \cref{tab:hyper-parameters_landscape}. The true parameters $(\beta^{\dag},\gamma^{\dag})$ are fixed as $(-21.6, 1.6)$, which are typical values for optical fibers.

\begin{table}[!htbp]
    \centering
    \begin{tabular}{c c c c c c c}
    \toprule
    symbol rate $1/T_{s}$ & $Z$ & $M$ & $N_{s}$ & zero-padding & sampling rate & \\
    \midrule
    100 GBaud & 80 km & 100 & 200 & 70/side & $2^6$ sps\\
    \bottomrule
    \end{tabular}
    \caption{Default setting for visualizing loss landscape. `sps' stands for samples per symbol.}
    \label{tab:hyper-parameters_landscape}
\end{table}

A few empirical observations can be made. Firstly, the landscape is non-convex and contains many local minimizers. Secondly, it exhibits anisotropicity, which is caused by different nature of the two parameters $\beta$ and $\gamma$. Thirdly, it has an unique global minimizer $(\beta^{*},\gamma^{*})$, which is the optimal solution of the optimization problem but may not be the true solution $(\beta^{\dag},\gamma^{\dag})$ of the inverse problem. Finally, the basin of attraction containing $(\beta^{*},\gamma^{*})$ is deeper and wider than those containing other local minimizers. Here we have borrowed the term \textit{basin of attraction} from theory of dynamical systems \cite{hirsch2012differential}. Each minimizer of the landscape can be regarded as an attractor of the gradient flow of the landscape. The basin of attraction containing any attractor is a set of points in the parameter space, which are drawn to the attractor by the gradient flow. For example, the deep blue elliptic region in the center of the landscape in \cref{fig:loss_landscape} is contained in a basin of attraction containing the global minimizer. 

\subsection{Effects of Model Architecture}
\label{ss:effect_model}
According to the empirical observation of the loss landscape, there is a unique global minimizer $(\beta^{*},\gamma^{*})$. If we have an `oracle' training algorithm, this minimizer is the optimal estimate we can find. On the other hand, this solution depends on model architecture and data, which are characterized by a set of hyper-parameters. Regardless of the specific algorithm for solving the inverse problem, we can get some knowledge about the relationship between the optimal estimate and these hyper-parameters, which is helpful for our understanding of the problem and also provides guidance for designing model and choosing training algorithms. 

\begin{figure}[!htbp]
    \centering
    \subfigure[number of layers]
      {\includegraphics[width=0.32\textwidth]{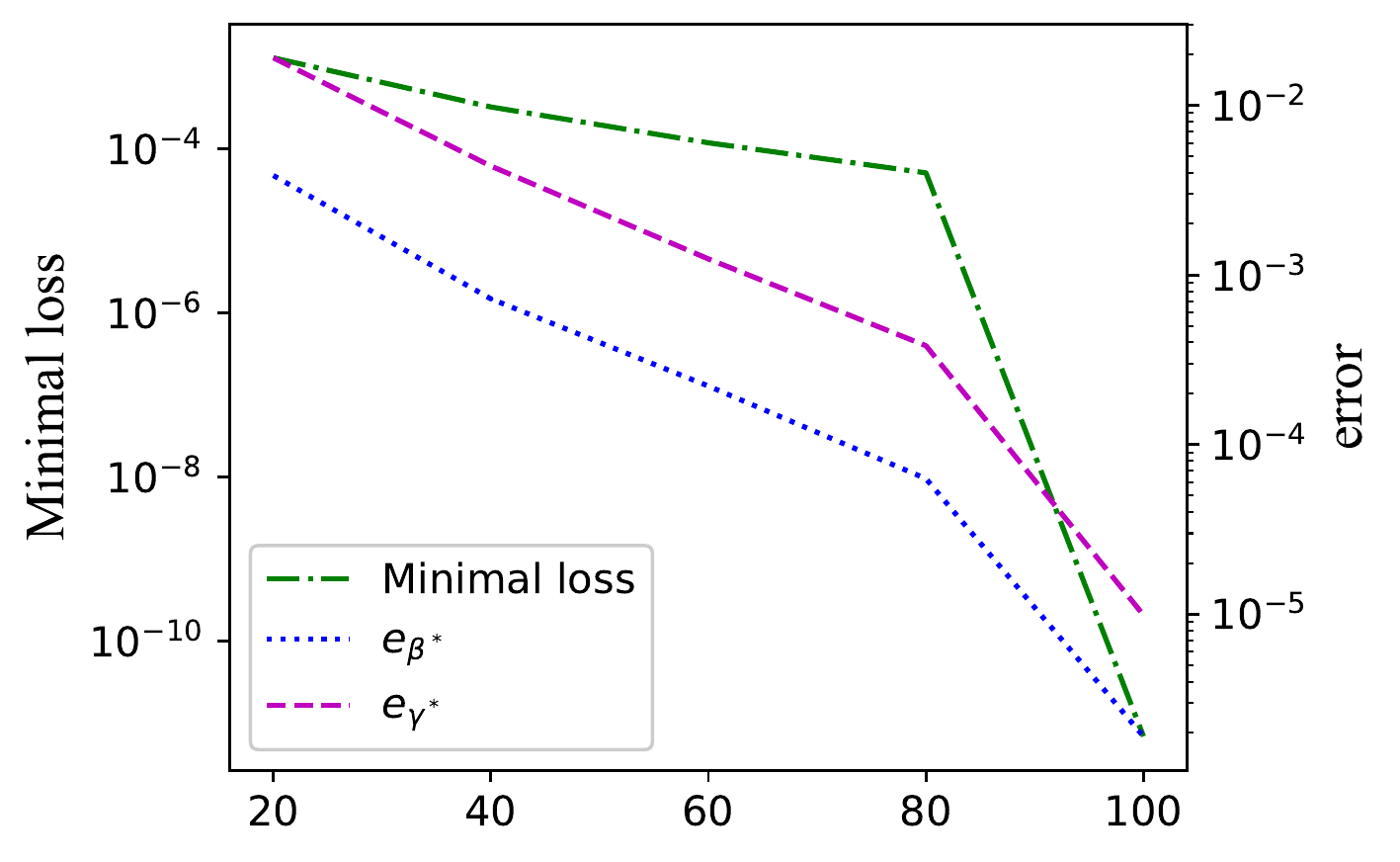}\label{fig:dependence_layers}}
    \subfigure[sampling rate]
      {\includegraphics[width=0.32\textwidth]{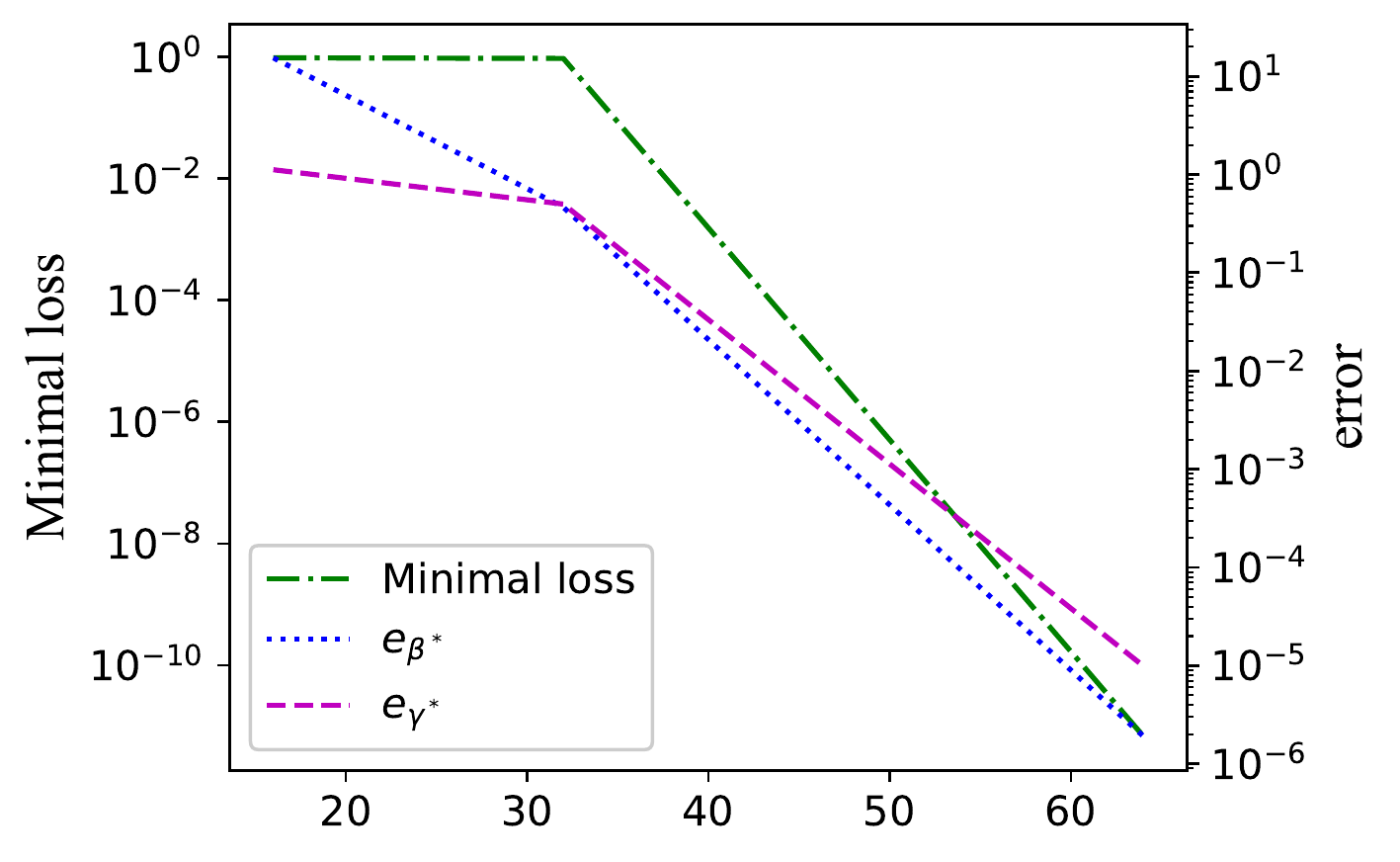}\label{fig:dependence_sampling_rate}}
    \subfigure[number of symbols]
      {\includegraphics[width=0.32\textwidth]{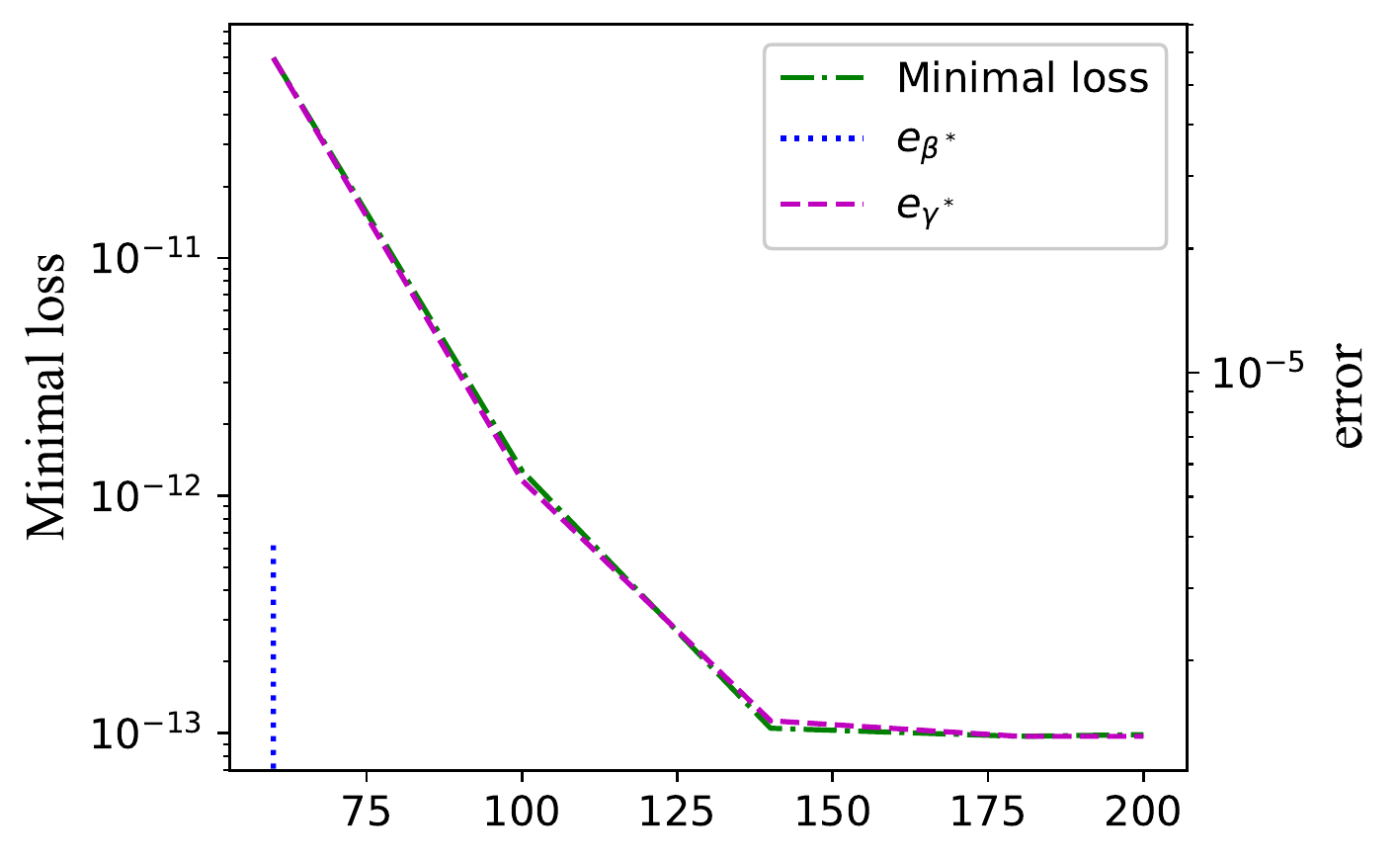}\label{fig:dependence_symbols}}
    \caption{Dependence of global minimum and minimizer on hyper-parameters of the model. The left axes indicate values of global minimum. The right axes indicate estimation errors of $\beta$ and $\gamma$ at global minimizer. The horizontal axes indicate values of hyper-parameters.}
    \label{fig:solution_dependence}
\end{figure}

Here we investigate the dependence of the optimal estimate on two hyper-parameters: the number of layers in NLS-Net and the sampling rate (sps), which characterize the depth and width of the NLS-Net respectively. The other settings are fixed as \cref{tab:hyper-parameters_landscape}. The results are shown in \cref{fig:solution_dependence}. The minimal loss $J(\beta^{*},\gamma^{*}):=\min_{\beta,\gamma} J(\beta,\gamma)$ and errors of the optimal estimate $(\beta^{*},\gamma^{*})$ are computed for each hyper-parameter respectively. The estimation errors of $(\beta,\gamma)$ are defined as $e_{\beta}=|\beta-\beta^{\dag}|$ and $e_{\gamma}=|\gamma-\gamma^{\dag}|$, where $\beta^{\dag}$ and $\gamma^{\dag}$ are ground truths simulated with the setting in \cref{tab:hyper-parameters_landscape}. 
The $(\beta^{*},\gamma^{*})$ in each setting is found by minimizing $J(\beta,\gamma)$ using gradient descent (GD) algorithm.

We trained NLS-Net for 20, 40, 60, 80 and 100 layers with GD until convergence. The results are summarized in \cref{fig:dependence_layers}. It tells us that both minimal loss and errors of optimal estimate decreases as the number of layers increases, and the decreasing rate increases. 
For sampling rate, we have chosen three different values: 16, 32 and 64 (samples per symbol). As we can see from \cref{fig:dependence_sampling_rate}, both minimal loss and estimation errors of optimal estimate decreases as sampling rate increases, and the decreasing rate increases. But 32 is too small since the optimal estimation errors are around $10^{-1}$. Hence 64 is the relatively proper choice.

From \cref{fig:dependence_layers} and \cref{fig:dependence_sampling_rate}, we can also have some empirical observations about the expressive power of NLS-Net. Recall that \cref{thm:convergence_rate} gives an explicit upper bound of the minimal loss, i.e. $\zeta^{2}+\tau^{s}$. But this bound does not converge to zero for fixed $\zeta$ or $\tau$, which is unsatisfactory in the view of machine learning. In contrast to the theoretical results, the numerical simulations summarized in \cref{fig:dependence_layers} and \cref{fig:dependence_sampling_rate} suggest that converging to zero is possible when one of the depth and width is fixed. In \cref{fig:dependence_layers}, the width is fixed as the depth increases. In \cref{fig:dependence_sampling_rate}, the depth is fixed as the width increases. In both results, we see the trend of converging to zero. Although we can not prove those convergence results rigorously, they provided support to use NLS-Net as a machine learning model for the problem.

Minimal loss is a measure of approximation error of output of NLS-Net. The uniqueness of global minimizer as observed in \cref{fig:loss_landscape} implies that as minimal loss converges to zero, the estimation errors $e_{\beta^{*}}$ and $e_{\gamma^{*}}$ should also converge to zero. Combined with the \cref{thm:convergence_rate}, it suggests strongly that the optimal estimation errors of $(\beta,\gamma)$ converge to zero as the depth and width of NLS-Net goes to infinity. This conclusion agrees with the numerical results in \cref{fig:dependence_layers} and \cref{fig:dependence_sampling_rate}. Therefore, the role of optimal approximation error (minimal loss) can exchange with that of optimal estimation error. 

We also did numerical experiments when both the depth and width of NLS-Net are fixed. The \cref{fig:dependence_symbols} shows the change of errors with respect to the number of symbols in the data. We can see that both minimal loss and optimal estimation errors of solution decreases as the number of symbols increases, but the decreasing rate decreases. They get saturated when the number of symbols exceed some critical value. It suggests that when we fix both depth and width of the NLS-Net, there is a nonzero lower bound of approximation error. Hence the expressive power of NLS-Net is limited by its depth and width.  Consequently, there is always a gap between the global minimizer of loss landscape and the true solution, no matter how much data we use.

\subsection{Effect of Data}
\label{ss:effect_data}
Solving the inverse problem is to find a good estimate of the truth parameter $(\beta^{\dag},\gamma^{\dag})$ of NLS-Net. Due to the limitation and uncertainty of collecting data, it is natural and necessary to consider the influence of data in the estimation of parameters. We need to answer the questions that how estimation error is affected by the length and randomness of data.

We recall that the input data \cref{eq:Ai} is generated from a sequence of random symbols $\{a_{k}\}\subset\cC$. Considering the randomness of data, the gap between the optimal estimation $(\beta^{*}, \gamma^{*})$ and the ground truth $(\beta^{\dag}, \gamma^{\dag})$ is better understand in terms of bias-variance decomposition. The bias is the gap between the ground truth and the mean of optimal estimations, which reflects the limitation imposed by the model. This part has been discussed in \cref{ss:effect_model}. The variance describes the concentration of $(\beta^{*}, \gamma^{*})$ around there mean, which is caused by the limitation of data.

\begin{figure}[!htbp]
    \centering
    \subfigure[$N_s=50$]
      {\includegraphics[width=0.45\textwidth]{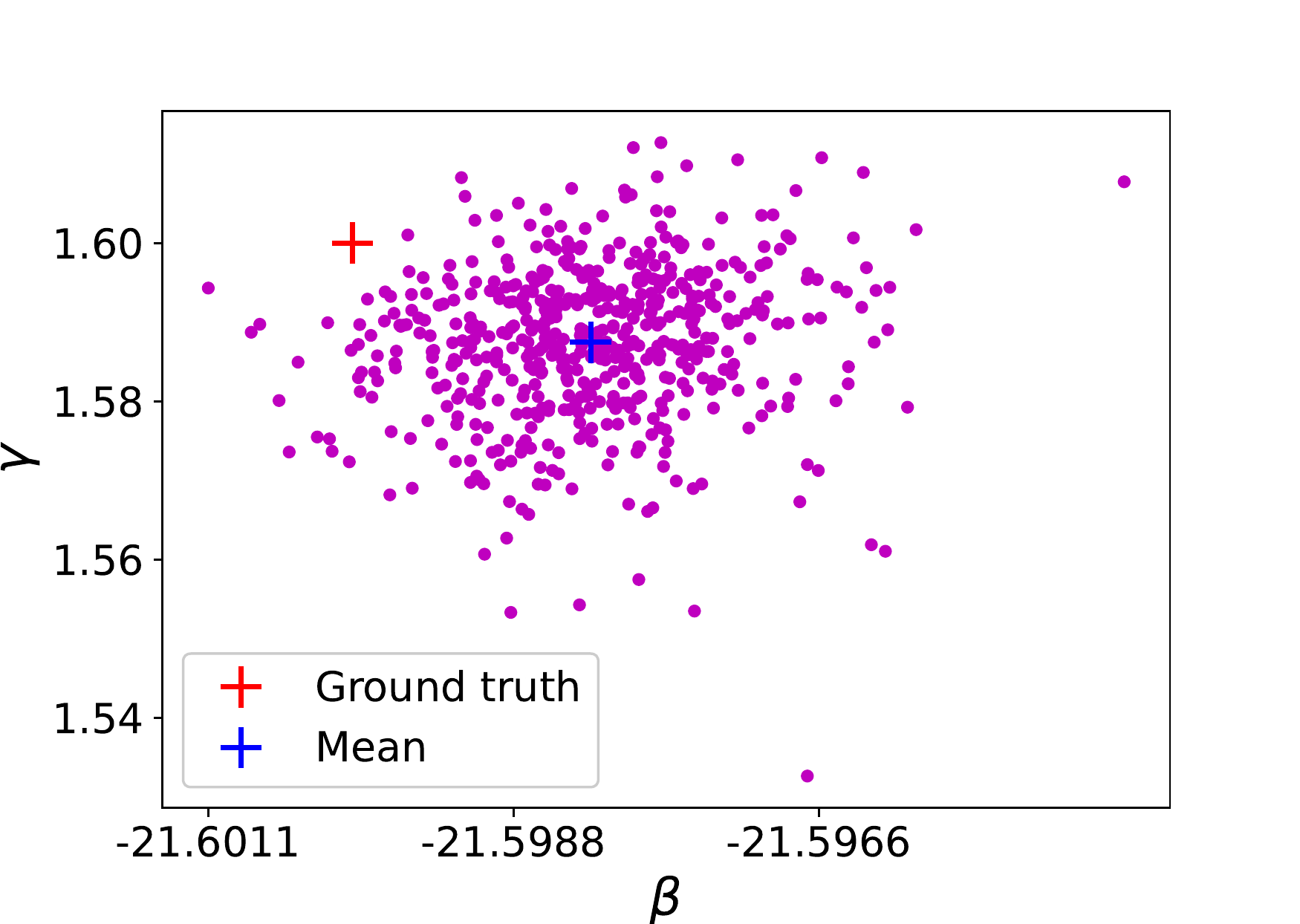}\label{fig:dis_50}}
    \subfigure[$N_s=100$]
      {\includegraphics[width=0.45\textwidth]{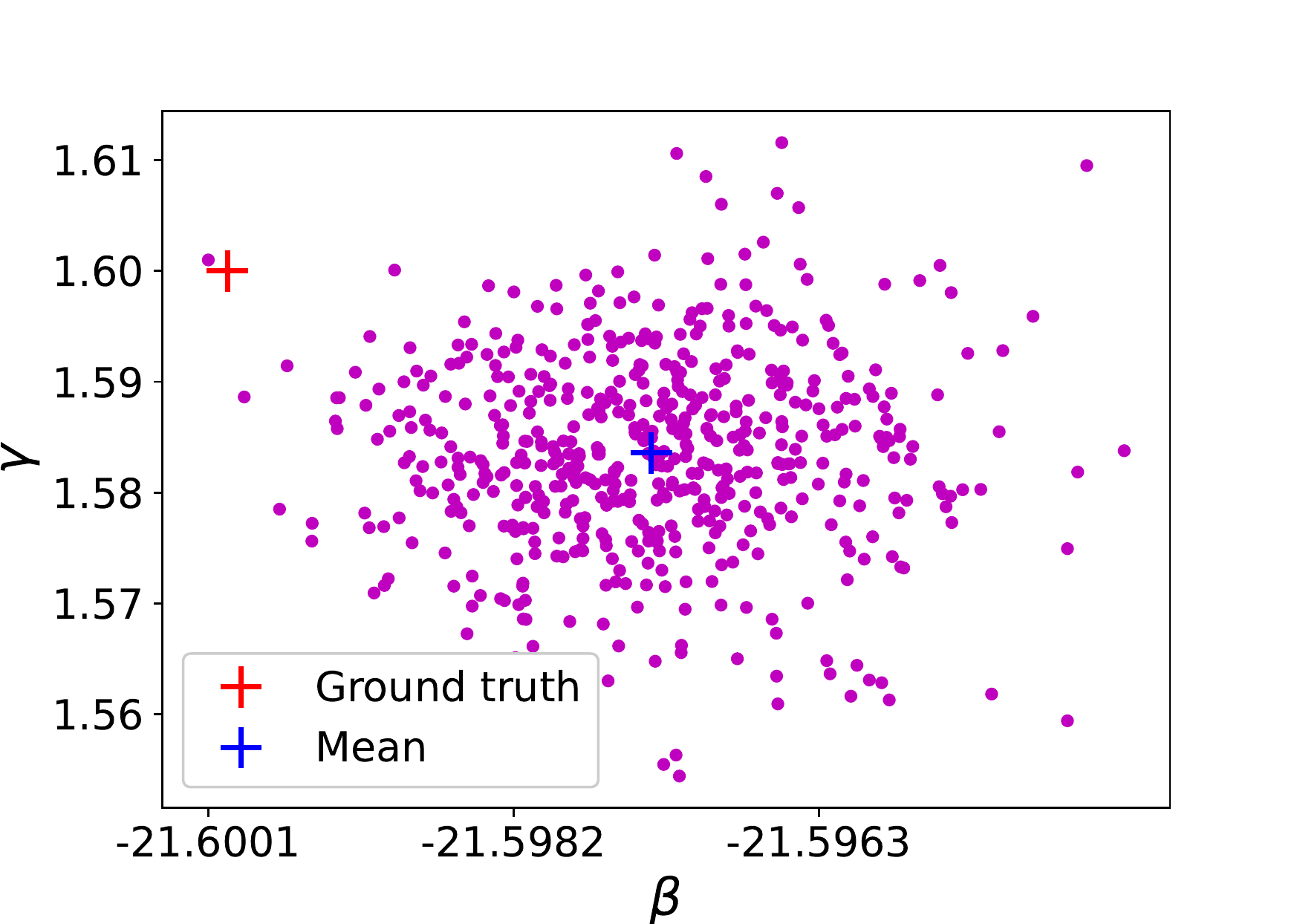}\label{fig:dis_100}}
    \subfigure[$N_s=150$]
      {\includegraphics[width=0.45\textwidth]{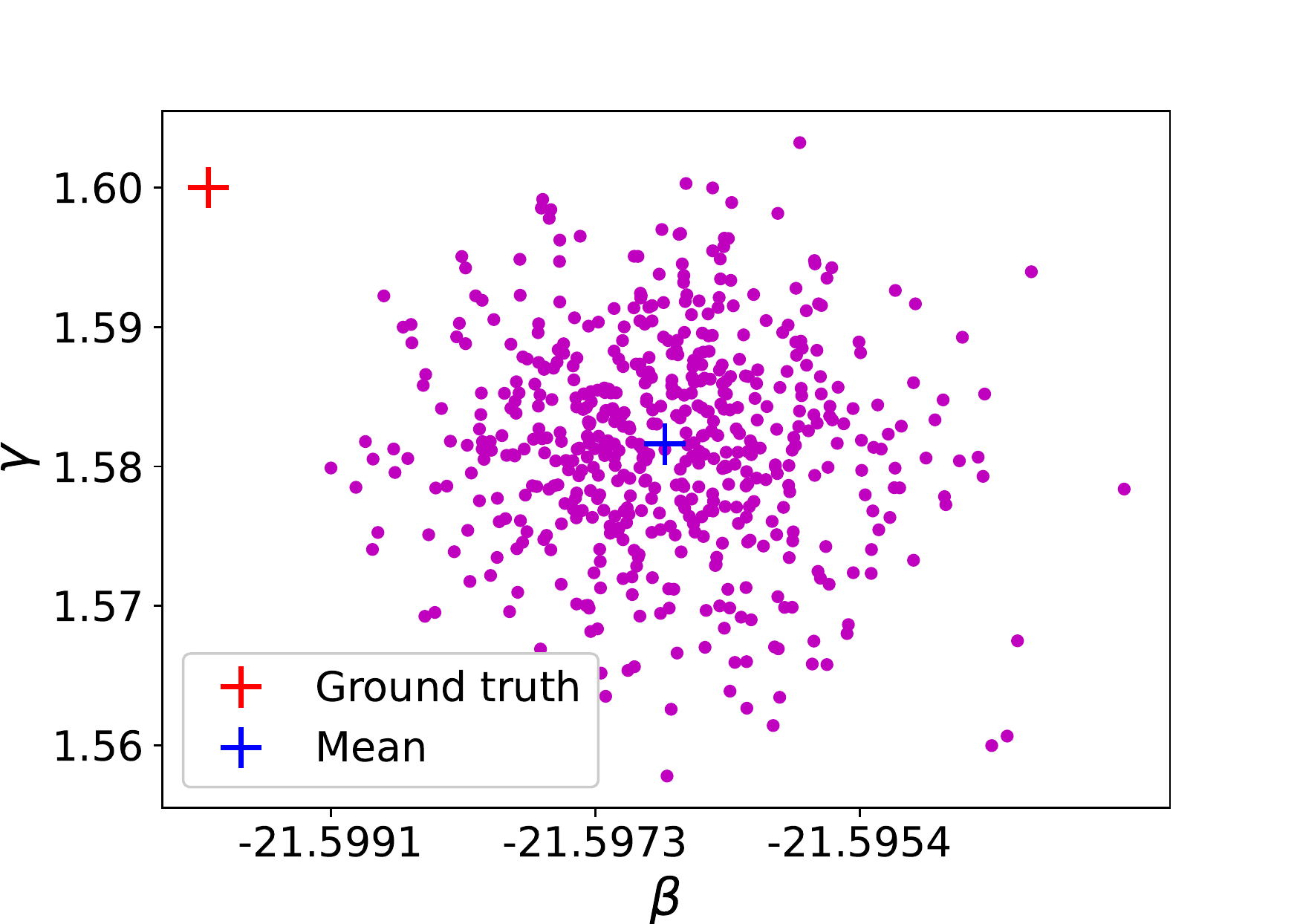}\label{fig:dis_150}}
    \subfigure[$N_s=200$]
      {\includegraphics[width=0.45\textwidth]{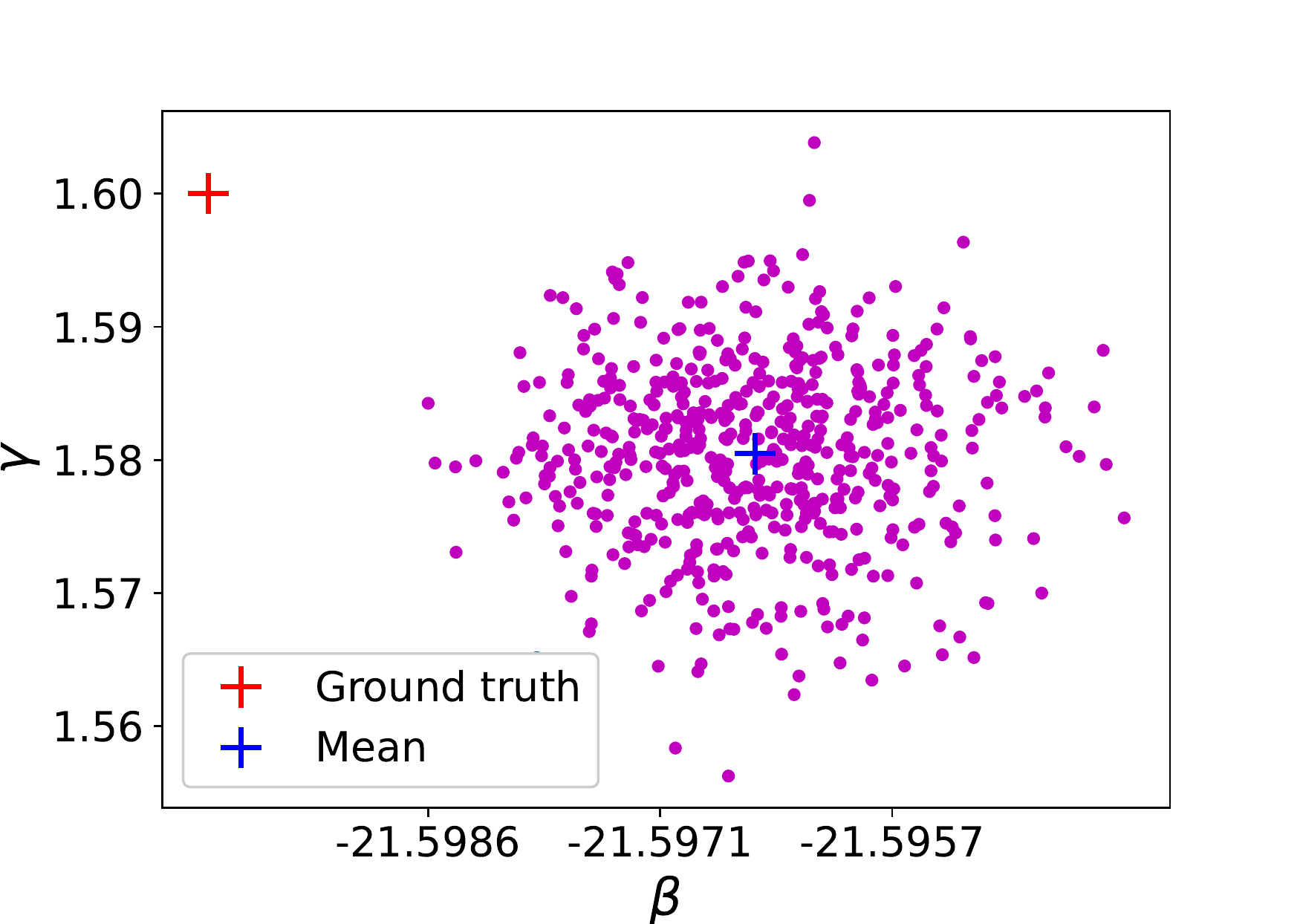}\label{fig:dis_200}}
    \caption{Distributions of global minimizers of random data, grouped by length of data.}
    \label{fig:minimizer_locations}
\end{figure}

\begin{figure}[!htbp]
    \centering
    \subfigure[Density of $\beta^{*}$]
      {\includegraphics[width=0.45\textwidth]{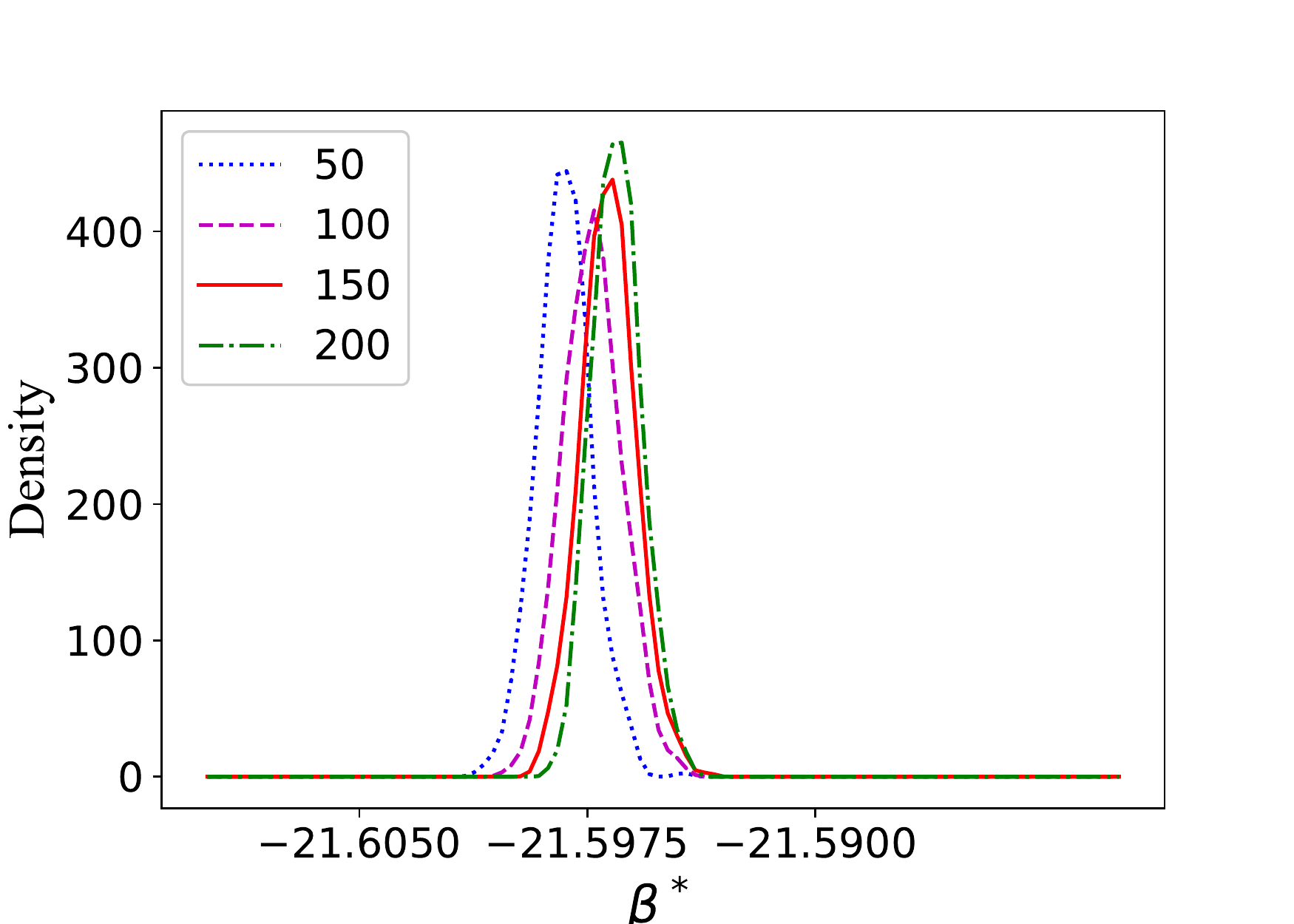}\label{fig:pdf-beta}}
    \subfigure[Density of $\gamma^{*}$]
      {\includegraphics[width=0.45\textwidth]{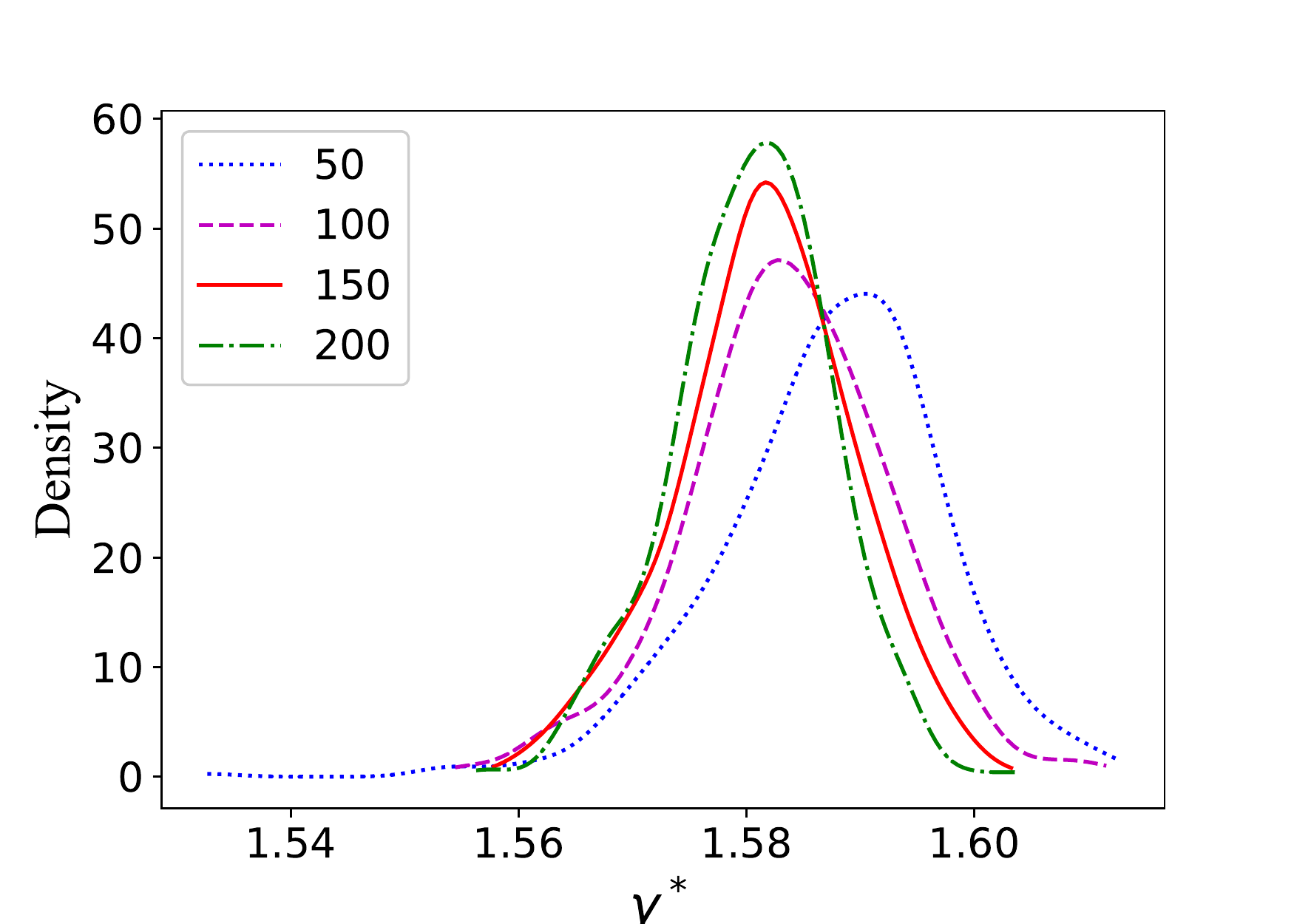}\label{fig:pdf-gamma}}
    \caption{Density functions of marginal distributions.}
    \label{fig:minimizer_density}
\end{figure}

We have drawn 4 groups of such sequences of random symbols, each contains 550 sequences of the same length. From the sequences we generate data under the same default setting \cref{tab:hyper-parameters_landscape}. Then we solve the inverse problem by training a 20-layer NLS-Net with each of the data. The optimal estimates corresponding to the data are plotted in \cref{fig:minimizer_locations}, grouped according to the length of data. The marginal density functions of $\beta^{*}$ and $\gamma^{*}$ are shown in \cref{fig:pdf-beta} and \cref{fig:pdf-gamma} respectively.
\begin{figure}[!htbp]
    \centering
    \subfigure[$\beta^{*}$]
      {\includegraphics[width=0.45\textwidth]{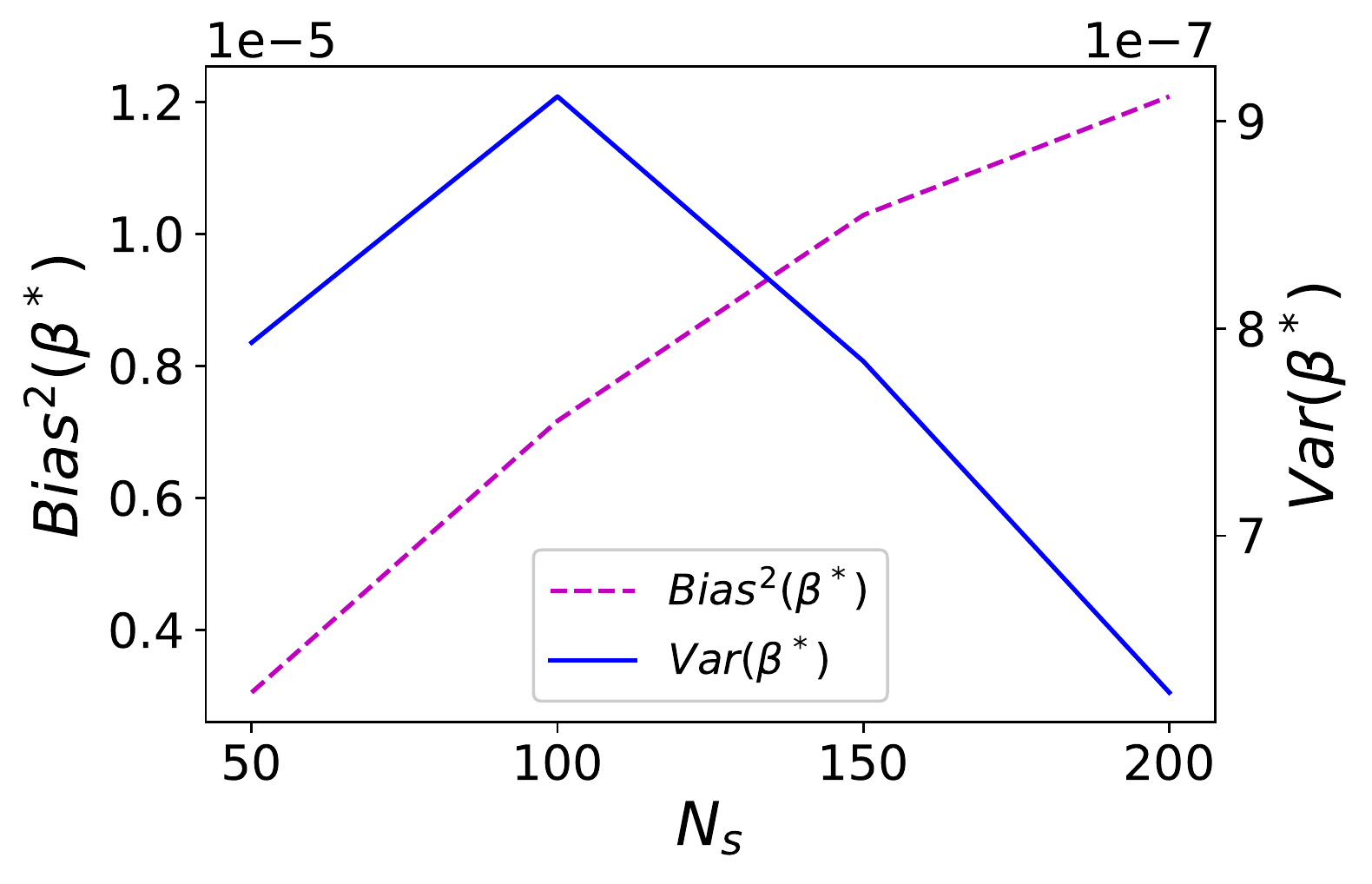}\label{fig:var-beta}}
    \subfigure[$\gamma^{*}$]
      {\includegraphics[width=0.45\textwidth]{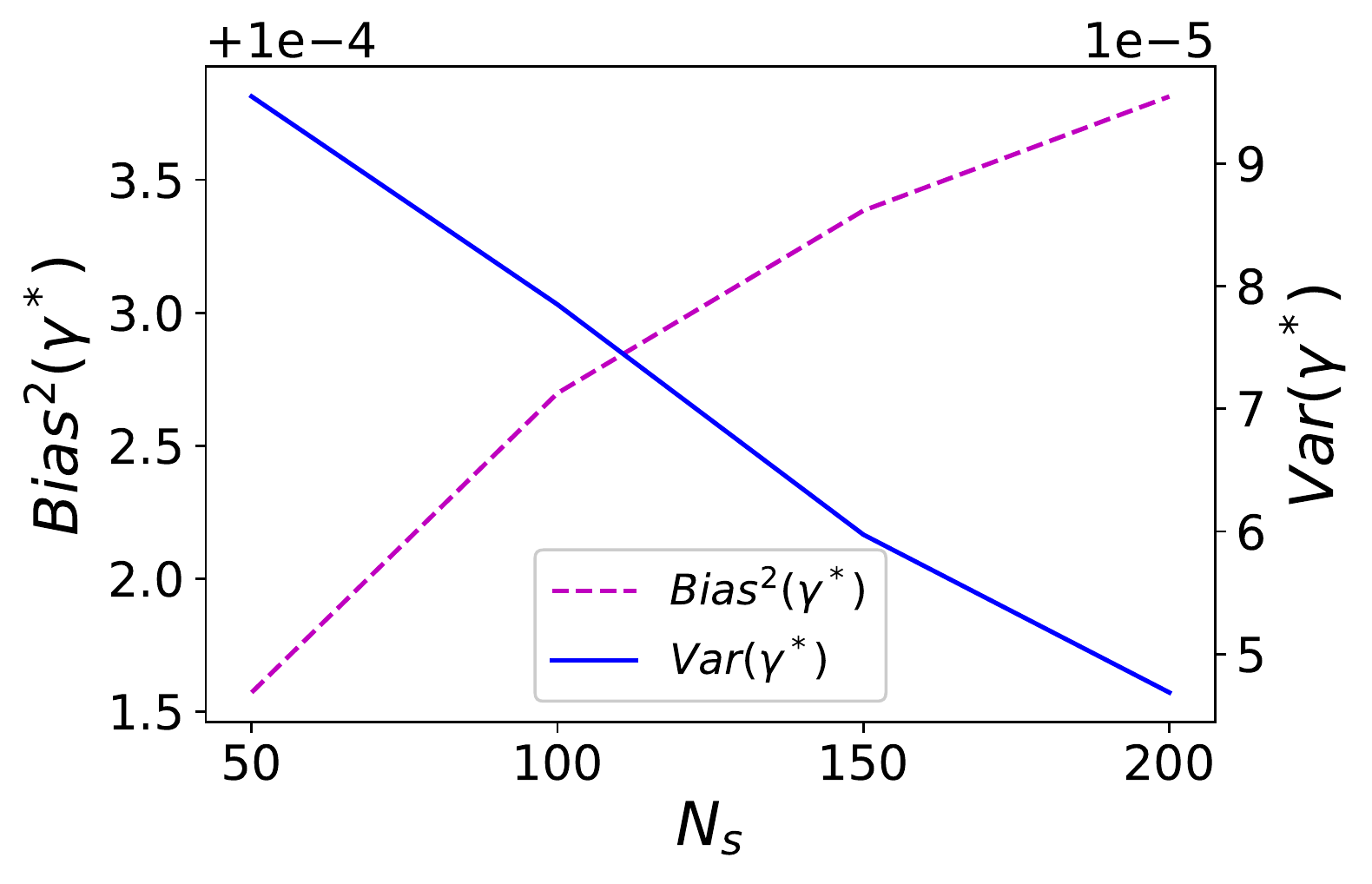}\label{fig:var-gamma}}
    \caption{Bias and variances of $\beta^*$ and $\gamma^*$.}
    \label{fig:minimizer_bias_var}
\end{figure}

\begin{table}[!htbp]
    \centering\small
    \begin{tabular}{c c c c c}
    \toprule
    $N_{s}$     & $50$ & $100$ & $150$ & $200$\\
    \midrule
    Mean of $\beta^*$    & $-21.5983$ & $-21.5973$ & $-21.5968$ & $-21.5965$\\
    \midrule
    Mean of $\gamma^*$  &  $1.5875$ & $1.5836$ &  $1.5816$ & $1.5805$\\
    \bottomrule
    \end{tabular}
    \caption{Empirical means of $(\beta^*,\gamma^*)$ corresponding to different numbers of symbols.}
    \label{tab:minimizer_mean}
\end{table}
\begin{table}[!htbp]
    \centering\small
    \begin{tabular}{c c c c c}
    \toprule
    $N_{s}$ & $50$ & $100$ & $150$ & $200$\\
    \midrule
    & $\begin{pmatrix} 0.0793& 0.122 \\ 0.122 & 9.56 \end{pmatrix}$
    & $\begin{pmatrix} 0.0913 & 0.0158 \\ 0.0158 & 7.87  \end{pmatrix}$
    & $\begin{pmatrix} 0.0786 & -0.0191 \\ -0.0191 & 5.99 \end{pmatrix}$
    &$\begin{pmatrix}  0.0625& 0.00383 \\ 0.00383 &  4.70  \end{pmatrix}$\\
    \bottomrule
\end{tabular}
    \caption{Covariance matrices corresponding to different numbers of symbols. Each matrix as an extra factor $10^5$.}
    \label{tab:minimizer_covariance}
\end{table}

We can see from \cref{fig:minimizer_locations} that the optimal estimates $(\beta^{*}, \gamma^{*})$ concentrate in a narrow region of parameter space. The value of empirical means of all groups are listed in \cref{tab:minimizer_mean}.
As we have already discussed in \cref{ss:effect_model}, there is always a gap between the ground truth and the optimal estimate, no matter how long data we use. These empirical observations are re-confirmed by \cref{fig:minimizer_locations}, in which there is a gap between the ground truth and the means of optimal estimates. This gap corresponds to the bias error in the bias-variance decomposition. See \cref{tab:minimizer_mean} for the values and see \cref{fig:minimizer_bias_var} for visualizations. These results shows that the bias increases as the length of data increases but it increases slower as more data is used. Hence it is predictable that the bias has a tendency to converge to the approximation error corresponding to the chosen $\zeta$.

The concentration of optimal estimates as shown in \cref{fig:minimizer_locations} tells us that most of these minimizers are closed to each other. This is quantitatively characterized by the empirical covariance matrix of $(\beta^{*}, \gamma^{*})$. See \cref{tab:minimizer_covariance} for the values and see \cref{fig:minimizer_bias_var} for visualizations.
These results show that for sufficiently long data, the variance of optimal estimates decreases as the number of symbols increases. That is, increasing data reduces variance. The covariance matrices also exhibit anisotropicity in the two parameters $\beta$ and $\gamma$. It indicates that changing data affects the estimation of $\beta$ less than estimation of $\gamma$. 
\section{Solving the Inverse Problem} \label{ss:solving}

Now we are to minimize the loss function \cref{eq:full_loss}. As shown in \cref{fig:loss_landscape}, it is a non-convex function.
Generally speaking, our task is easier than the training of DNNs, since the optimization problem involves only two variables $(\beta, \gamma)$. On the other hand,  our task has an extra difficulty. In most tasks of DNN-training, one only needs to minimize the loss function with \textit{any} solution, while in our task one is required to find the \textit{unique} solution. Even if a solution is equally good with the true solution in generating the output data, it may not be the true solution.

Nevertheless, the scope of this paper is restricted to convex optimization. The development of algorithms for the general non-convex optimization of $J$ is left as future work. In this section, we use prior knowledge to restrict the optimization variables $(\beta,\gamma)$ to a small region and compare the performance of several training algorithms that are popular in deep learning. More specifically, we choose the starting point of $(\beta,\gamma)$ within the basin of attraction containing the ground truth $(\beta^{\dag},\gamma^{\dag})=(-21.6,1.6)$. Hence the optimization problem is essentially convex. The compared algorithms are the gradient descent (GD) with momentum, the Adam \cite{Kingma_2014}, the Adadelta \cite{adadelta} and RMSprop \cite{RMSprop}.

All these algorithms are gradient-based. But we didn't use the stochastic version of these algorithms here. The gradient of the loss function with respect to optimization variables can be computed through the well-known backpropagation (BP) procedure, which is nothing but the application of chain rule on neural networks. It is most clearly illustrated through computation graph. See \cref{fig:computation_graph}.
With the help of this computation graph, we can easily apply chain rule and compute the gradients of $J$ with respect to $\beta$ and $\gamma$. The name \textit{backpropagation} means that the gradient is computed backwardly along all paths connecting the output variable with $\beta$ (or $\gamma$). Thanks to deep learning software such as PyTorch \cite{pytorch}, the computation of gradient can be done automatically.

The numerical experiments are done with three groups of data: noiseless, noisy and denoised. Since we have studied the loss landscape of the problem with noiseless data, it is natural to solve the problem in the same setting to complete the study. Then we add observation noise to the data and solve the problem again to examine the performance of the method and algorithms in more practical setting. Finally, we process the noisy data with simple denoising method and solve the problem again. We find that denoising is effective to improve the quality of solution.

\begin{figure}[!htbp]
    \centering
      {\includegraphics[width=0.32\textwidth]{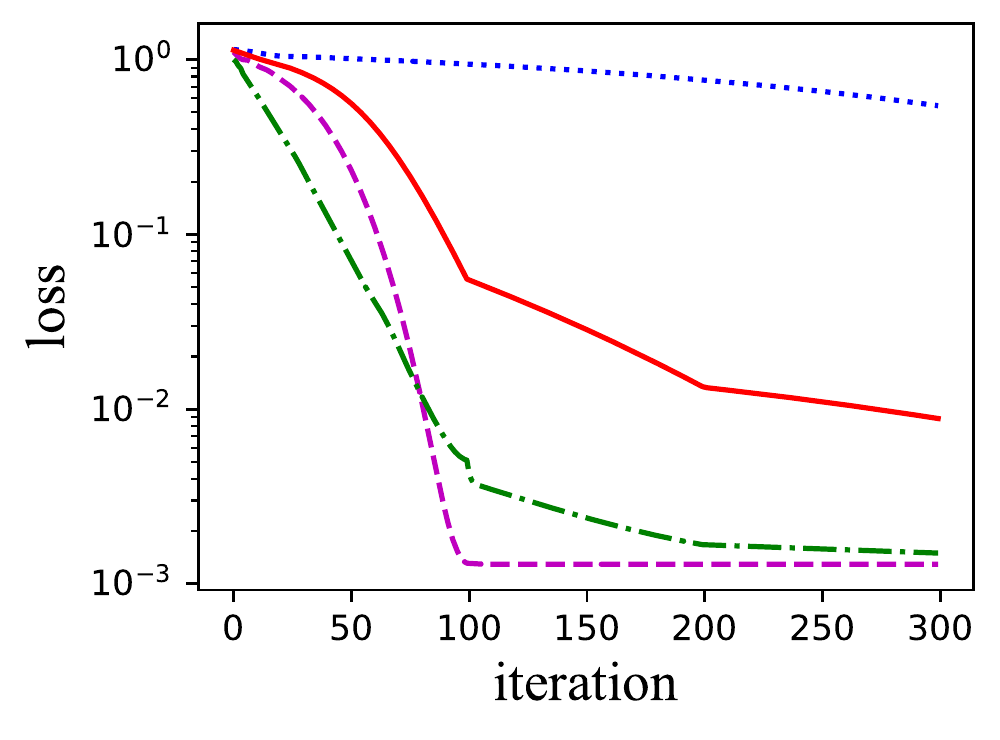}\label{fig:loss_pure}}
      {\includegraphics[width=0.32\textwidth]{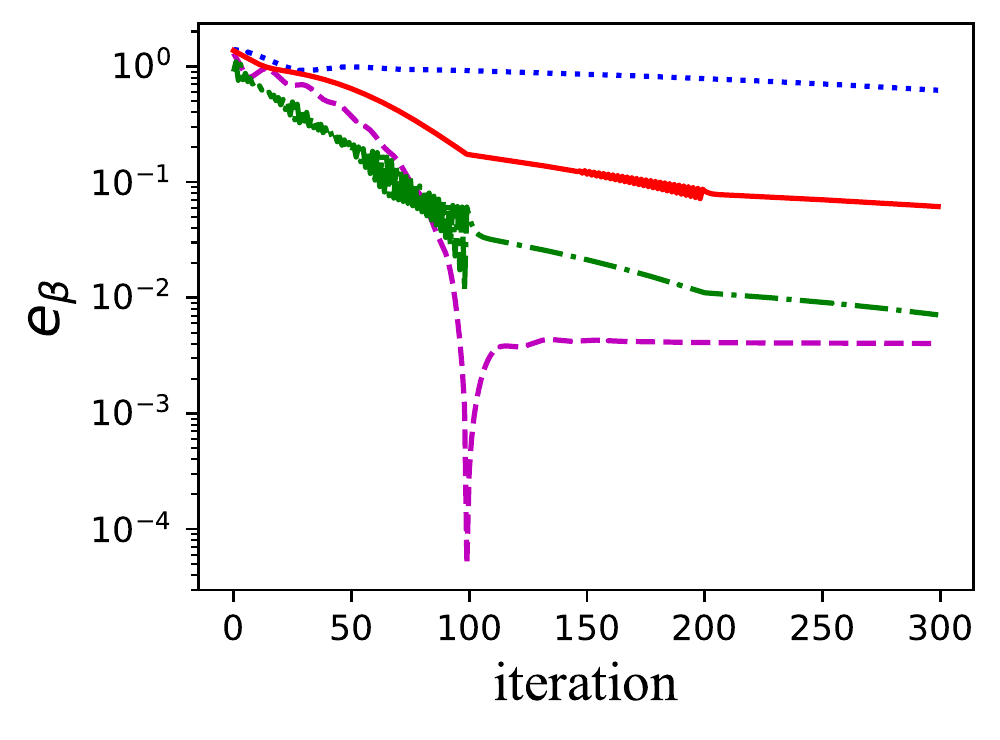}\label{fig:beta_pure}}
      {\includegraphics[width=0.32\textwidth]{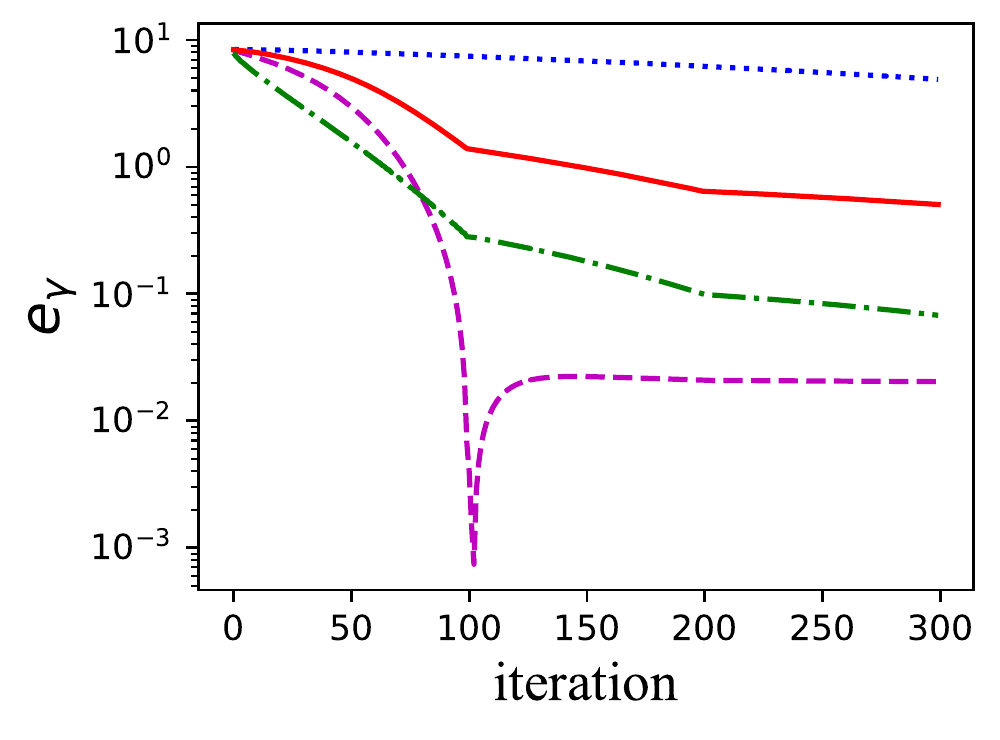}\label{fig:gamma_pure}}
      {\includegraphics[width=0.32\textwidth]{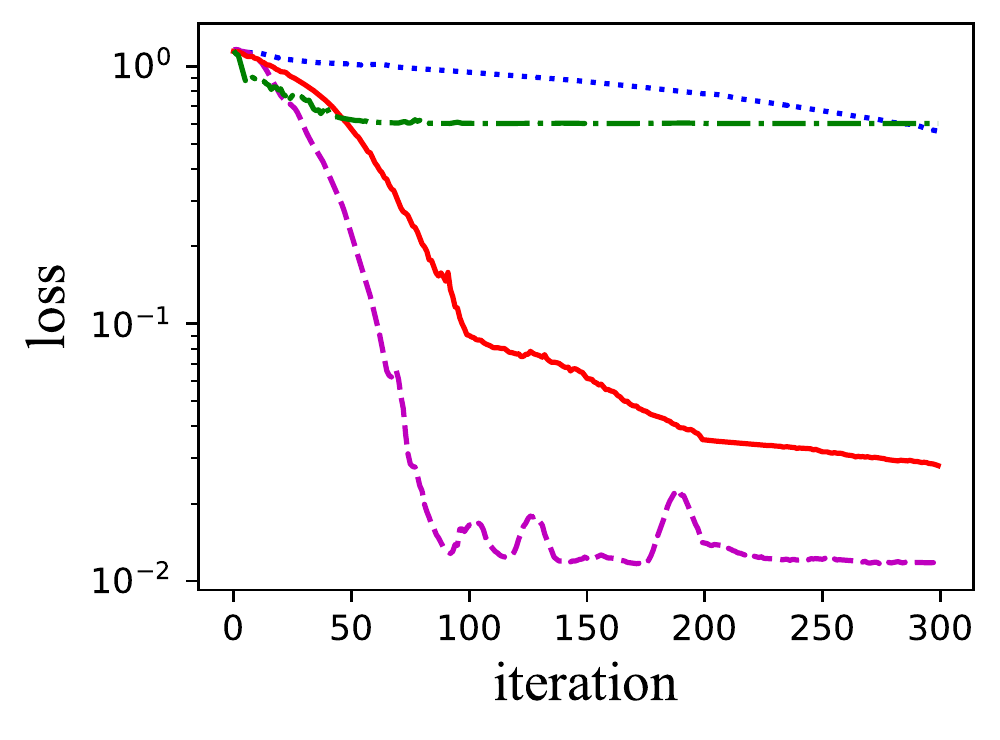}\label{fig:loss_noisy}}
      {\includegraphics[width=0.32\textwidth]{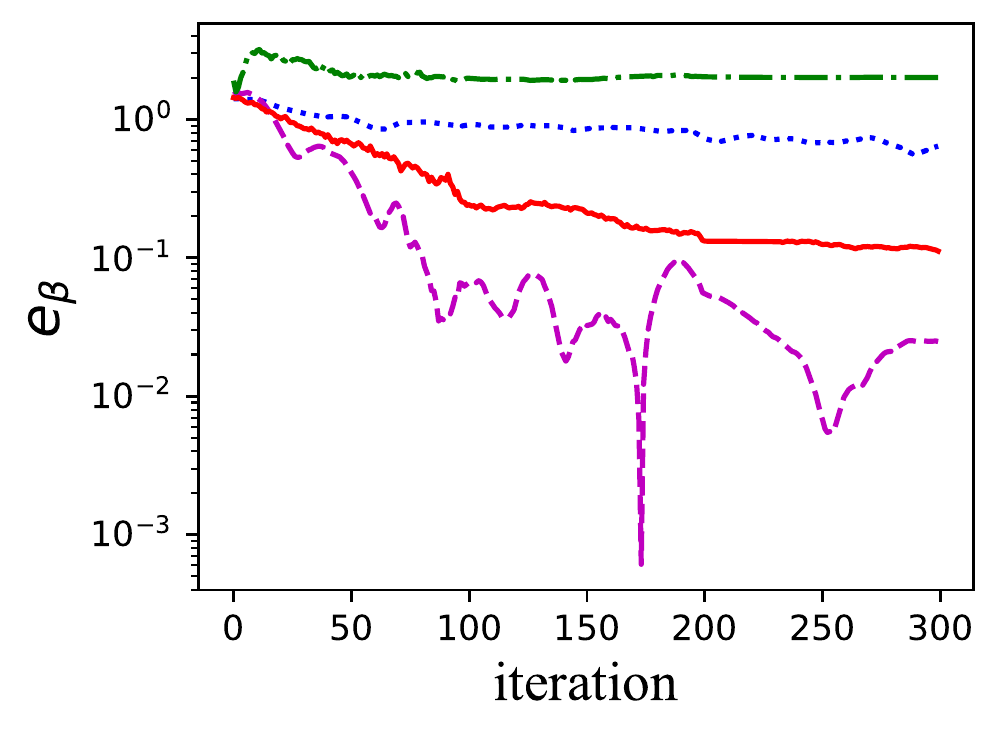}\label{fig:beta_noisy}}
      {\includegraphics[width=0.32\textwidth]{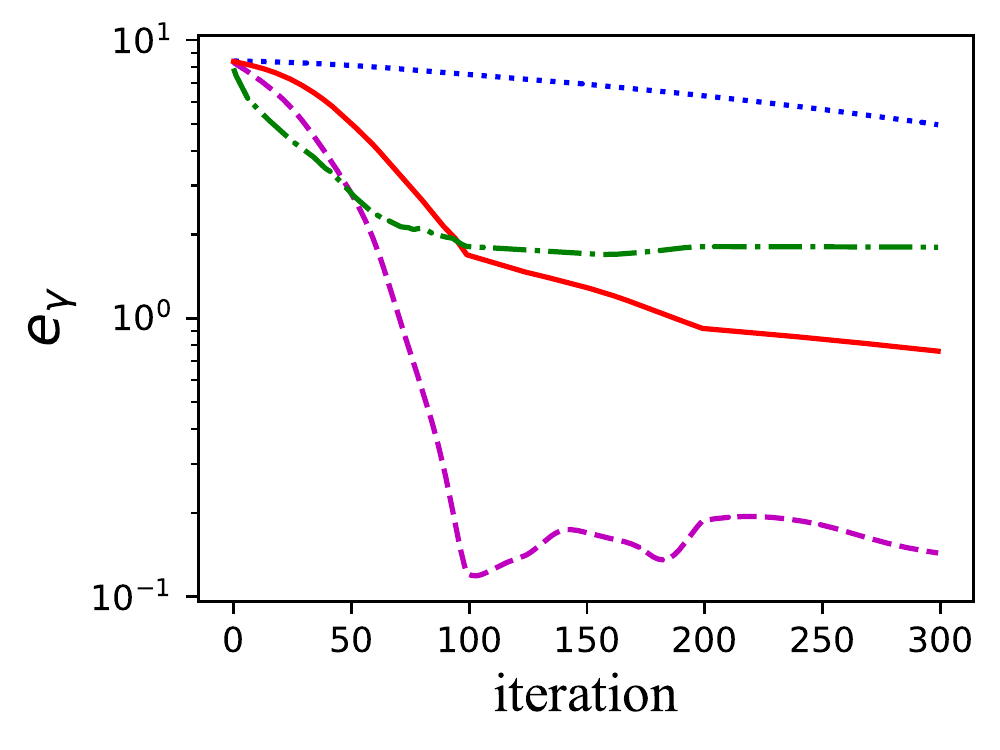}\label{fig:gamma_noisy}}
      {\includegraphics[width=0.32\textwidth]{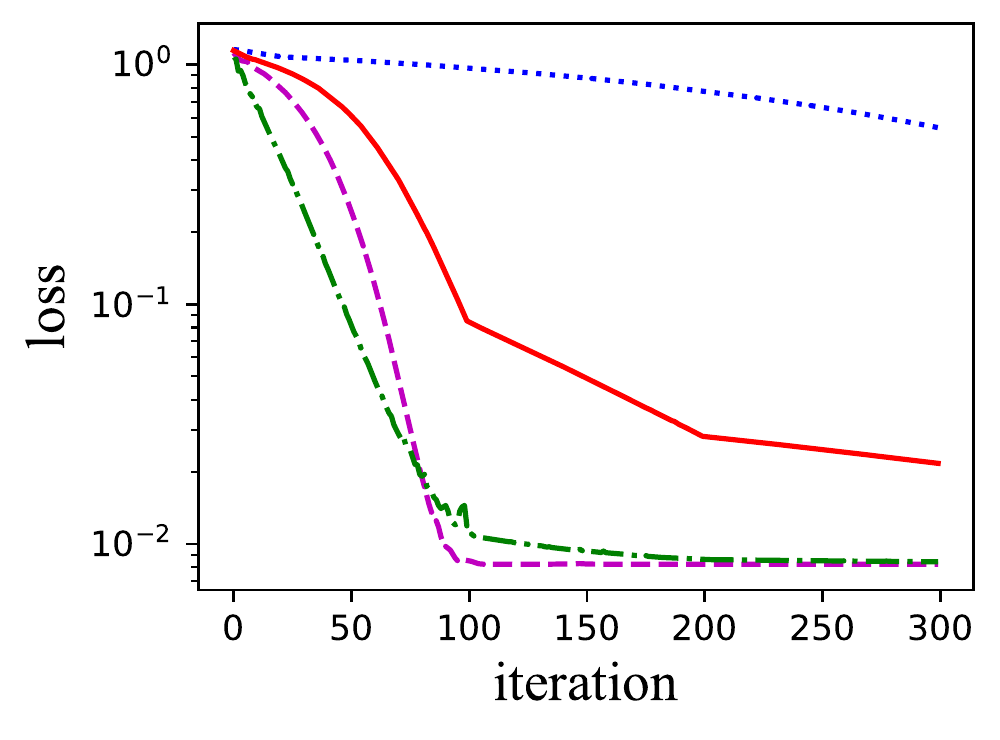}\label{fig:loss_filter}}
      {\includegraphics[width=0.32\textwidth]{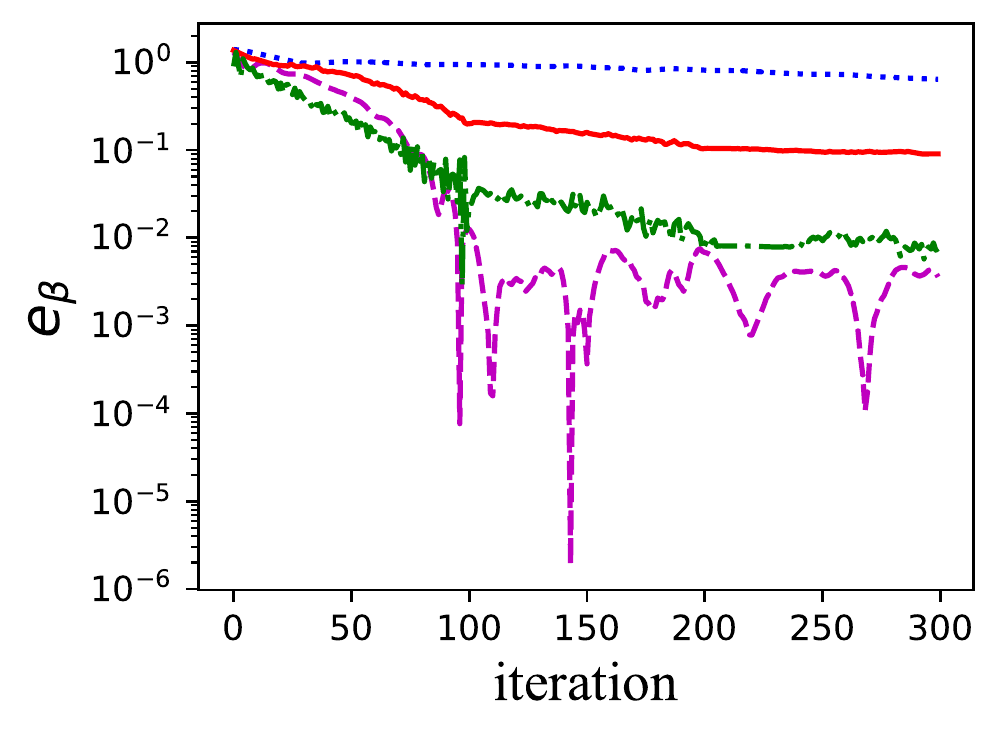}\label{fig:beta_filter}}
      {\includegraphics[width=0.32\textwidth]{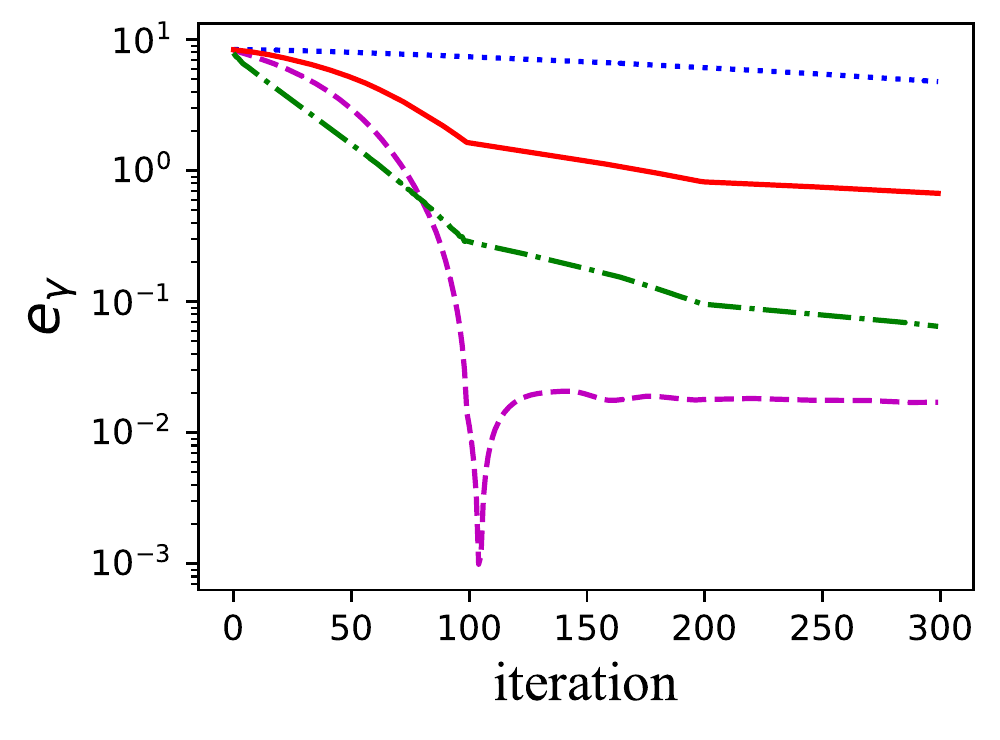}\label{fig:gamma_filter}}
      {\includegraphics[width=0.5\textwidth]{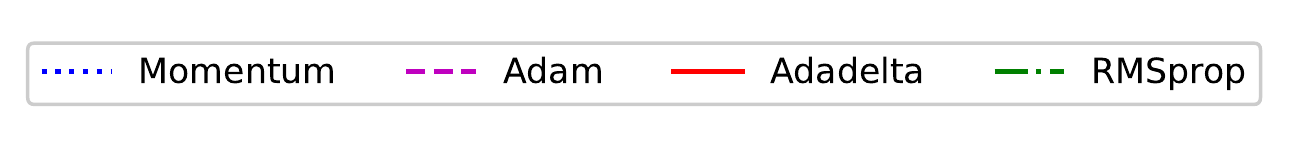}}
    \caption{Comparison of different algorithms. First row: with noiseless data; second row: with noisy data; third row: with denoised data. First column: loss function; second column: error of $\beta$; third column: error of $\gamma$.}
    \label{fig:training_history_comparison}
\end{figure}

The data are generated in the same default setting as described in \cref{tab:hyper-parameters_landscape}. The $h$ is defined in \cref{eq:pulse_shape} with $\rho=0.1$ and $P=1$. Noise is added according to \cref{ss:data} with SNR=200. The denoising procedure is to apply the $h$ as a filter to both $\tilAi$ and $\tilAo$, i.e., to convolve $\tilAi$ and $\tilAo$ with $h$. The starting point of $(\beta,\gamma)$ is set to $(-23,10)$ in all experiments.
The training history of our model on three groups of data are shown in \cref{fig:training_history_comparison}.

Several observations are made about the results shown in \cref{fig:training_history_comparison}.  Let's look at the first column, which shows the history of loss during training. On all the three groups of data, the convergence of Adam is fastest, then it is Adadelta, and that of GD with momemtum is the slowest. Their order is independent of noise. In contrast, RMSprop is affected significantly by noise. The estimation errors of parameters, showed in the second and third columns of \cref{fig:training_history_comparison}, exhibit similar order.
Secondly, the comparison of the three rows shows that the converged loss and estimation errors on the noisy data (middel row) are higher than those on the noiseless and denoised data, while the latter two are comparable. It suggests that noise is harmful and denoising is helpful. Hence the denoising procedure is quite effective.

\begin{figure}[!ht]
    \centering
    \subfigure[Outputs of training data.]{
	\includegraphics[width=0.8\textwidth]{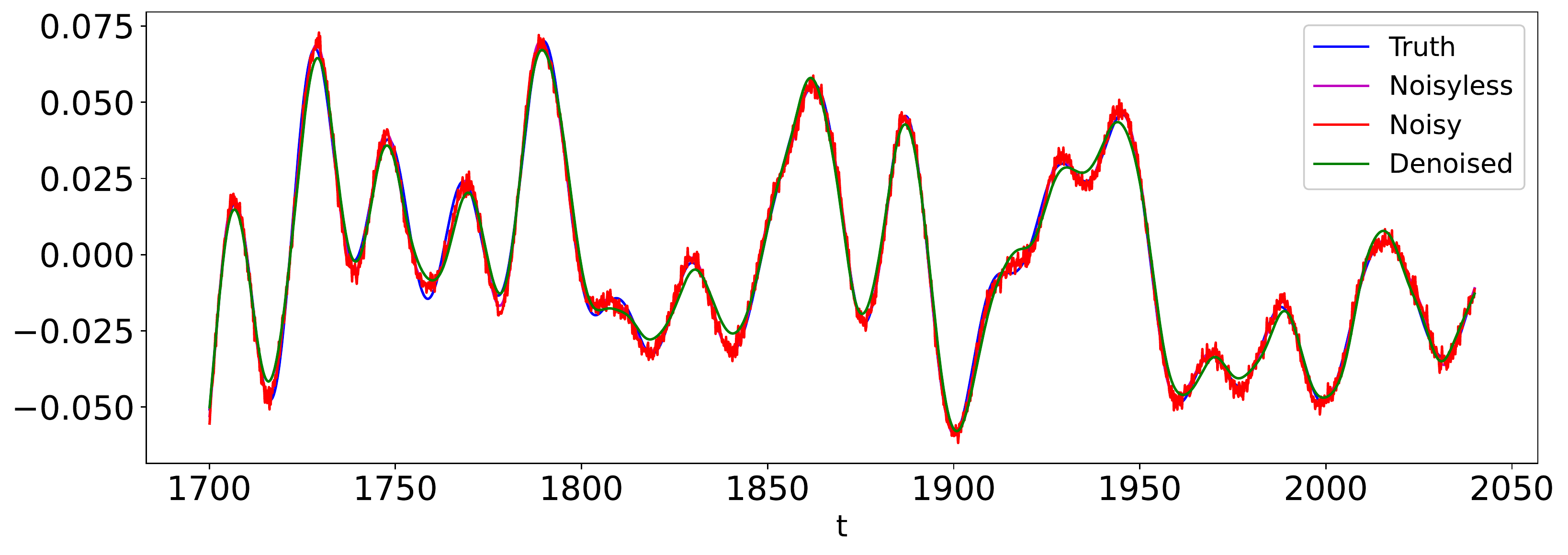}}
	\subfigure[Outputs of testing data.]{
	\includegraphics[width=0.8\textwidth]{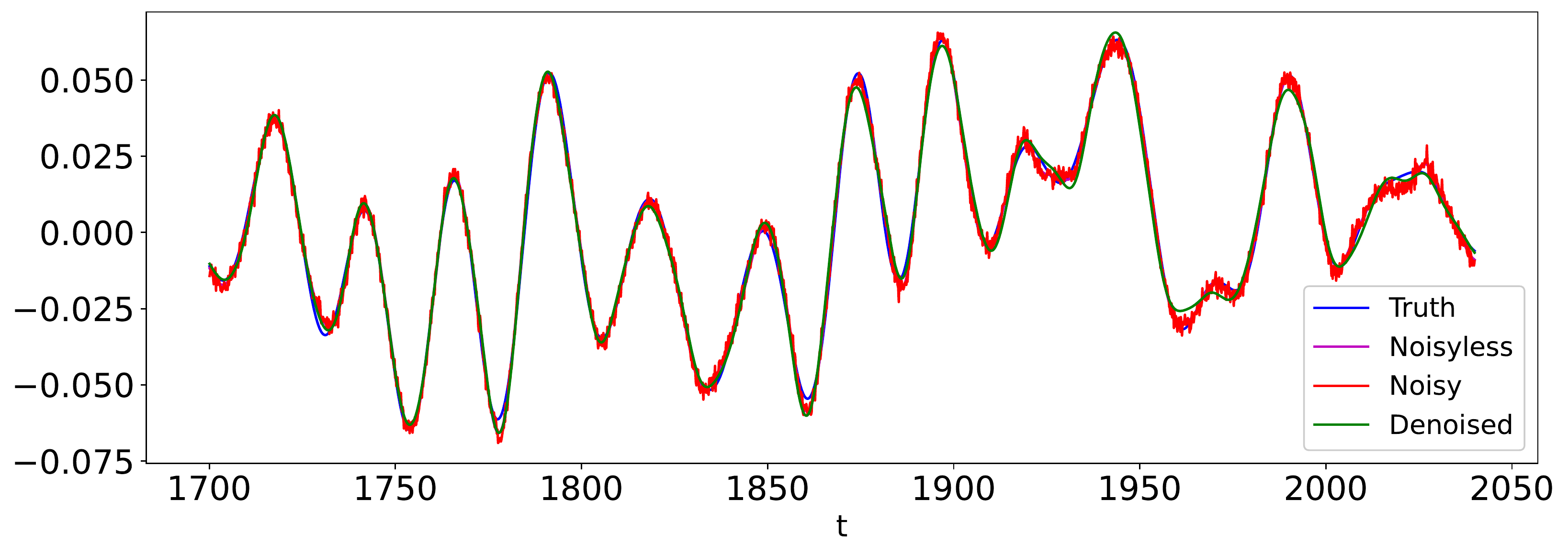}}
	\caption{Outputs of a trained NLS-Net. Only the real part of each data is shown.}
    \label{fig:comparison_outputs}
\end{figure}

To get a better sense of how good are the solutions, we can visualize a representative part of the output of the trained NLS-Net corresponding to the final state of  \cref{fig:training_history_comparison}. See \cref{fig:comparison_outputs} for an illustration of the outputs corresponding to the training and testing data respectively. The training data is used to estimate $(\beta^{\dag},\gamma^{\dag})$ by Adam. The testing input data is in the same form as \cref{eq:Ai_comp}, which is feed into the trained NLS-Net to obtain its output. Numerical experiments are did on noiseless, noisy and denoised data. The results are compared with the true outputs. We can see that on all data the outputs are very close to the true outputs. This is an indirect evidence that the solution of the inverse problem is very accurate.

\section{Discussions} \label{ss:discussion}

In this paper, we formulate and solve an inverse problem of NLSE as a learning problem of the NLS-Net. In this paper, however, we focus on the estimation of parameters and didn't discuss the generalization problem, i.e., prediction on new input data. That will be studied in our future works.

The key of our method is to consider the SSFM solver of NLSE as a re-parameterized convolutional neural network. The NLS-Net is inherently explainable and its expressive power is guaranteed by the theory of SSFM, which gives an upper bound of the approximation error of NLS-Net, but this bound does not converge to zero for NLS-Net with a fixed depth or width. However, our empirical study suggests that the actual approximation error of such NLS-Net goes to zero if one of the depth and width is fixed. Numerical results also show that there is a nonzero lower bound of approximation error when both of depth and width are fixed. 
According to our empirical observations, the loss landscape of the learning problem of NLS-Net has an unique global minimizer. It implies that the convergence of minimal loss is equivalent to the convergence of estimation error of parameters. Hence the above conclusions about the approximation error can be transferred to estimation error.

We also studied the dependence of the optimal estimate on different data. The randomness of data causes sampling error for estimation of parameters. 
The optimal estimates obtained from different data concentrate within a region of the parameter space around their mean, and they get close to the mean as the amount of data increases. At the same time, they get close to each other. However, increasing the amount of data does not reduce the gap between the mean estimate and the ground truth, which is limited by the expressive power of the NLS-Net. Hence we have observed the bias-variance decomposition in the estimation of parameters. 

Finally, we compare several training algorithms on noiseless, noisy and denoised data. It is found that simple denoising procedure such as filtering works quite well. It can almost cancel the effect of observation noise. In addition, numerical results suggest that the Adam algorithm is the best among others. Visualizations shows that the obtained estimate is very accurate.

\appendix
\section{NLSE with Nonzero Fiber Attenuation}
In the main part of this paper, we studied the inverse problem of NLSE with chromatic dispersion and Kerr nonlinearity, where the fiber attenuation effect is absent. 
Generally, we can consider the NLSE \cref{eq:NLSE_IVP2} in which there is another coefficient $\alpha$ describing the fiber attenuation. We show that under some general conditions, the coefficient $\alpha$ can be obtained in closed form without solving optimization problem. The generalized NLSE is
\begin{align}
\begin{cases}
\label{eq:NLSE_IVP2}
    \frac{\p A(t,z)}{\p z}=- \frac{\alpha}{2} A(t,z)-\frac{i\beta}{2} \frac{\p^{2} A(t,z)}{\p t^{2}}+i\gamma|A(t,z)|^{2} A(t,z), & t\in \bbR, z\in [0,Z],\\
    A(t,0)=\Ai(t), & t\in \bbR,
\end{cases}
\end{align}
where $\Ai$ is the initial signal at the transmitter with $z=0$. Similarly as in \cref{ss:abs_formulation}, $\Ao(t)$ is defined to be the data at the receiver with $z=Z$. We have the following
\begin{theorem}
In the IVP \cref{eq:NLSE_IVP2}, suppose that $\Ai\in H^{1}(\bbR\to\bbC)$ and let $A\in H^{1}_{t}C_{z}(\bbR\times[0,Z]\to\bbC)$ be the solution of the IVP. Denote $\Ao(t)=A(t,Z)$ for all $t$, then we have
\begin{align}
\label{eq:estimation_alpha}
    \alpha = \frac{2}{Z}\log\frac{\|\Ai\|_{L^{2}}}{\|\Ao\|_{L^{2}}}.
\end{align}
\end{theorem}
\begin{proof}
Denote the complex conjugate of $A$ as $A^{*}$. Using \cref{eq:NLSE_IVP2} it is easy to verify that
\begin{align}
\label{eq:energy_ineq}
    \frac{\p}{\p z}\int_{\bbR}|A(t,z)|^{2}dt
    =-\alpha \int_{\bbR}|A(t,z)|^{2}dt
    +\beta \imag \left(A^{*}(t,z)\frac{\p A(t,z)}{\p t}\right)_{t=-\infty}^{\infty}.
\end{align}
Since $A\in H^{1}_{t}C_{z}(\bbR\times[0,Z]\to\bbC)$, for any fixed $z\in[0,Z]$, there is
\begin{align}
    \int_{\bbR} \left|\imag \left(A^{*}(t,z)\frac{\p A(t,z)}{\p t}\right)\right| dt \le \left(\int_{\bbR} |A(t,z)|^{2} dt \int_{\bbR} \left|\frac{\p A(t,z)}{\p t}\right|^2 dt\right)^{1/2}<\infty.
\end{align}
Hence the last term in \cref{eq:energy_ineq} vanishes. Integrate \cref{eq:energy_ineq} in $z$, we get
\begin{align}
    \|\Ao\|_{L^{2}}^{2} = \|\Ai\|_{L^{2}}^{2} e^{-\alpha Z}.
\end{align}
It is easy to verify that \cref{eq:estimation_alpha} holds.
\end{proof}

\section*{Acknowledgements}
This research did not receive any specific grant from funding agencies in the public, commercial, or not-for-profit sectors.
\bibliographystyle{plain}
\bibliography{references}
\end{document}